\documentclass[12pt]{article}
\usepackage[centertags]{amsmath}
\usepackage{amsfonts}
\usepackage{amssymb}
\usepackage{amsthm}
\usepackage{latexsym}
\usepackage{amsmath}
\usepackage{amsfonts}
\usepackage[active]{srcltx}
\usepackage{comment}
\usepackage{indentfirst}

\usepackage{graphicx}

\DeclareGraphicsExtensions{.eps,.pdf,.jpg,.png}
\usepackage{epstopdf}
\usepackage{subfig}
\usepackage{xcolor}
\usepackage{ulem}

\numberwithin{equation}{section}

\newtheorem{example}{Example}[section]
\newtheorem{theorem}{Theorem}[section]

\newtheorem{lemma}{Lemma}[section]

\newtheorem{remark}{Remark}[section]
\newcommand{\eproof}{{\mbox{\ }~\hfill
\mbox{\large $\Box$} \par \vskip 10pt}}

\newcommand{\pf}{\noindent{\bf Proof}}

\begin{document}
\title{
Clustered eigenvalue problem for glassy state relaxation and its inverse problem
}

\author{Shuli Chen\thanks{%School of Mathematics, Nanjing 210096, P. R. China 
(Email: sli\_chen@126.com).}
\and
Maarten V. de Hoop \thanks{Simons Chair in Computational and Applied
Mathematics and Earth Science, Rice University, Houston TX, USA (Email: mdehoop@rice.edu).}
\and
Youjun Deng\thanks{School of Mathematics and Statistics, Central South University,Changsha 410083, P.R. China (youjundeng@csu.edu.cn).}
\and
Ching-Lung Lin\thanks{Department of Mathematics, National Cheng-Kung University, Tainan 701, Taiwan (Email:
cllin2@mail.ncku.edu.tw).}
\and Gen Nakamura\thanks{Department of
Mathematics, Hokkaido University, Sapporo 060-0808, Japan and Research Center of Mathematics for Social Creativity, Research Institute for Electronic Science, Hokkaido University, Sapporo 060-0812, Japan (Email: gnaka@math.sci.hokudai.ac.jp).}}

\maketitle

\begin{abstract}
\noindent For computational convenience, a Prony series approximation of the stretched exponential relaxation function of homogeneous glasses has been proposed (\cite{MM}), which is the extended Burgers model known for viscoelasticity equations.  The authors of \cite{LS} initiated a spectral analysis of glass relaxation along this line, and gave some numerical results on clusters of eigenvalues. A theoretical justification of the results and development of further numerical studies were left open. In this paper, we provide a complete theoretical justification of their results and their numerical verification. Besides these, we solve an inverse spectral problem for clusters of eigenvalues associated with the glass relaxation.

\bigskip


\end{abstract}

\renewcommand{\theequation}{\thesection.\arabic{equation}}

\section{Introduction}

Glassy states can be formed in a variety of materials with different types of bonding properties such as oxides, alloys, molecules, and polymers. 
% In condensed matter physics, the nature of glasses and the glass transition are very important subjects to examine. 
Due to the non-equilibrium nature, atoms in glasses exhibit collective motions, notably viscoelasticity under deformation. This behavior is relevant for understanding the properties of glasses (see, for example, \cite{Song}). 

Non-exponential relaxation is a common property of glassy materials. For homogeneous glasses, the following empirical non-exponential function called stretched exponential function'', $G(t)$ say, is used. It is given as
\begin{equation}\label{s_ exp_func}
G(t)=\text{exp}\left\{-\left(\frac{t}{\tau}\right)^\beta\right\},  
\end{equation}
where $t\ge 0$ is time, $\tau>0$ is the relaxation time and $0<\beta<1$ is the stretching exponent. It is mathematically convenient to approximately replace $G(t)$ by a Prony series, which is of the form,
\begin{equation}\label{Prony}
G(t)=\sum_{i=1}^N s_i e^{-r_i t},
\end{equation}    
where $r_i>0,\,s_i>0$ for $1\le i\le N$ with $N\in{\mathbb N}$. Without loss of generality, we assume that
\begin{equation}\label{r's order}
0<r_1<r_2< \cdots <r_N.    
\end{equation}
It is reported in \cite{MM} that, away from $t=0$, the Prony series with a sufficient number of terms can accurately approximate the stretched exponential function. 

In this paper, we consider the initial boundary value problem associated with the viscoacoustic wave equation, with the above mentioned relaxation function, in space dimension one, given as
\begin{equation}\label{IBP}
\left\{
\begin{array}{rcl}
\partial_t^2 u(t, x) &=& \displaystyle{\int_0^t G(t-\tau)\partial_\tau  \partial_x^2 u(\tau, x)\,d\tau}\ \text{in $(0,\infty)\times\Omega$}, \\
u(t,0)&=&0,\quad \partial_x u(t,\frac{\pi}{2})=0\ \text{in $(0,\infty)$},\\[0.25cm]
u(0,x)&=&0\,\,\text{in $\Omega$},
\end{array}
\right.
\end{equation}
with $\Omega = (0,\pi/2)$. Here,
% $u$ denotes the displacement of the beam, and 
we have assumed that the density of the material is constant, equal to $1$, and that $G$ is independent of $x$ for simplicity. Integrating the right-hand side of the equation of motion in \eqref{IBP} by parts, we get
\begin{equation}\label{motion1}
\partial_t^2 u(t,x)=G(0)\partial_x^2 u(t,x)+\int_0^t G'(t-\tau)\partial_x^2 u(\tau,x)\,d\tau,    
\end{equation}
where we used the notation, $G'(t):=\partial_t G(t)$. Substituting the expression for the relaxation function, yields
\begin{equation}\label{motion2}
\partial_t^2 u(t,x)=D\,\partial_x^2 u(t,x)-\sum_{i=1}^N b_i\int_0^t e^{-r_i(t-\tau)}\partial_x^2 u(\tau,x)\,d\tau,    
\end{equation}
where $D=\sum_{i=1}^N s_i$ and $\,b_i=s_i\,r_i,\,\,1\le i\le N$. We omit specifying where the equations are satisfied unless it is unclear from the context. 

There are two main theories that provide the properties of relaxation. One is the (parametric) theory of spring-dashpot models (see \cite{DKLNT} and references therein), and the other is the theory of fading memory (see \cite{OR}, \cite{SJ} and references therein). Equation \eqref{motion2} is nothing but the equation of motion for one of the well-known spring-dashpot models, called the extended Burgers model (abbreviated by EBM). To guarantee that the solutions of \eqref{IBP} decay exponentially as $t\rightarrow\infty$, we need to connect a spring in parallel to this model. In that case, we have $D>h$ with $h:=\sum_{i=1}^N b_i/r_i$ (see \cite{DKLNT, DLN}). However, here, we will consider all cases,
\begin{equation*} % \label{cond_D}
D>h,\,\,D=h\ \text{and}\ D<h.
\end{equation*} 

\begin{comment}
\section{An approximate model for our problem}
Let $h>0$, $b_i>0$ for $1\leq i \leq N$,
and $0<r_1<r_2< \cdots <r_N$ such that
\begin{equation}\label{n1.2}
\begin{aligned}
\sum_{i=1}^N \frac{b_i}{r_i}=h.
\end{aligned}
\end{equation}
Consider the following Boltzmann-type viscoelastic equation asscociated with the extended Burgers model (abbreviated by EBM) {\color{blue} --- couldn't we do this in higher dimensions on a cube? we could discuss tis in a remark}
\begin{equation}\label{n1.1}
\begin{aligned}
u_{tt}=&D\partial_x^2 u - \sum_{i=1}^N b_i \int_0^t e^{-r_i(t-s)} \partial_x^2 u\, ds
\end{aligned}
\end{equation}
with the boundary condition
\begin{equation}\label{bc}
u(t,0)=0,\,\,\partial_x u(t,\frac{\pi}{2})=0,
\end{equation}
where $u=u(t,x),\,t\ge0,\,x\in\overline\Omega$, $u_{tt}=\partial_t^2 u$. Here, $\overline\Omega$ is the closure of $\Omega$ which is $[0,\pi/2]$ in our case.
\end{comment}
We introduce the shorthand notation $\xi u:=\partial_x^2 u $ and $I(r):=\int_0^t e^{-r(t-\tau)} \partial_x^2 u(\tau,x) \, d\tau $, to formulate an eigenvalue problem. With this notation, \eqref{motion2} takes the form,
\begin{equation}\label{n1.3}
\begin{aligned}
u''-D\xi u +  \sum_{i=1}^N b_iI(r_i)=0.
\end{aligned}
\end{equation}
% where $u'':=(u')'$ and “$'$" denotes the $t$-derivative.
We note that the following relation between $\xi$ and $I(r)$ holds,
\begin{equation}\label{n1.4}
\begin{aligned}
I'(r)=\xi u-r I(r).
\end{aligned}
\end{equation}
Then, introducing new dependent variables, $v=u'$ and $w_i=I(r_i)$ for $1\leq i\leq N$, we can write \eqref{n1.3} in the form of a system of partial differential equations that are of first order in time; we refer to this system, given by
\begin{equation}\label{n1.5}
\left\{
\begin{array}{rcl}
u'&=&v,\\
v'&=&D\xi u +  \sum_{i=1}^N b_iw_i,\\
w_i'&=&\xi u -r_iw_i, \qquad\quad\ 1\leq i \leq N,
\end{array}
\right.
\end{equation}
as the augmented system. The matrix-vector representation of this system is given by
\begin{equation}\label{n1.6}
\begin{aligned}
U'=A_{N+2}U,
\end{aligned}
\end{equation}
where $A_{N+2}$ is a $N+2$ square matrix of the second-order partial differential operators in $x$ and $U=(u,v,w_1,\cdots,w_n)^{\mathfrak{t}}$ with $\mathfrak{t}$ denoting the transpose.

Consider $A_{N+2}$ as a densely defined closed operator on the $N+2$ number of products of $L^2(\Omega)$, denoted by $L^2(\Omega)^{N+2}$, with the domain
\begin{equation}\label{domain of A}
D(A_{N+2}):=\left\{
\begin{array}{ll}
U=(u,v,w_1,\cdots,w_N)^{\mathfrak{t}}\in H^2(\Omega)\times H^1(\Omega)\times L^2(\Omega)^N:\\
\qquad\qquad\qquad u(t,0)=\partial_x u(t,\pi/2)=0,\,t>0
\end{array}
\right\},
\end{equation}
where $H^s(\Omega),\,s=1,2$ are the $L^2(\Omega)$-based  Sobolev spaces of order $s=1,2$. Then, it is natural to analyze the spectrum of $A_{N+2}$, specifically, the eigenvalue problem:
\begin{equation}\label{n1.7}
\begin{aligned}
\lambda U=A_{N+2}\,U.
\end{aligned}
\end{equation}

A comprehensive numerical study on giving an approximate correspondence between the stretched exponential function and its approximation by the Prony series is given in \cite{MM}; here, it is noted that spectral information of $A_{N+2}$ with $D=h$ can be used for estimation of this function from data. General data analysis for the estimation of the stretched exponential function can be found in \cite{Phillips, Song} and the references therein, where the data can be spectral data obtained by the dynamical mechanical analysis (abbreviated by DMA) instruments.  
% There is a lot of data analysis on the appropriate choice of $\beta,\,\tau$ of the stretched exponential function (see \cite{Phillips}, \cite{Song}, and the references therein).  However, it is very hard to find a paper on the eigenvalue problem \eqref{n1.7} other than \cite{LS}, which considered the case $D=h$. Needless to say that spectral information of $A_{N+2}$ with $D=h$ can be used for a data analysis of the stretched exponential function via the work on \cite{MM}. For instance, the spectral information can be obtained as measured data by using the DMA (dynamic mechanical analyzer).

Due to the difficulty that comes from $A_{N+2}$ being non-self-adjoint, the authors of \cite{LS} used the method given in the next section to obtain not necessarily all the eigenvalues but clusters of these. The authors provided compelling numerical results but without proof. This papers aims to give a comprehensive analysis with proofs with applies to the EBM.
% and left open problems to provide more numerical results and theoretical proof on them. This paper aims to solve these problems completely, and provide a basis to do the data analysis for the stretched exponential function describing the relaxation of glassy materials by using the mentioned approximate correspondence.
Although we consider a bounded interval in our analysis, we note that it can be applied to the case $\Omega \subset {\mathbb R}^d$ with $d > 1$ being a bounded domain with a Lipschitz smooth boundary, and $\partial_x^2$ is replaced by a second-order positive operator defined in $L^2(\Omega)$. This is because the operator has an eigenfunction expansion and its eigenvalues $\{\lambda_\ell\}_{\ell=1}^\infty$ increasingly labeled, counting multiplicity and satisfying $\lambda_\ell \sim c_0\,\ell^{d/2}$ for $\ell \gg 1$, where $c_0>0$ is a constant independent of $\ell$ \cite{Agr, CH}. 

The remainder of this paper is organized as follows. Section 2, we introduce the reduced eigenvalue problem with a parameter $k\in{\mathbb N}$ while describing the method mentioned above. The method can transform \eqref{n1.7} into an eigenvalue problem for a matrix $A_{N+2}^k$ which gives clustered eigenvalues of $A_{N+2}$. Then, in Section 3, we give the characterization of roots for the characteristic equation of $A_{N+2}^k$ and its limiting counterpart as $k\rightarrow\infty$. In Sections 4, 5, we present the speed of convergence of eigenvalues of $A_{N+2}^k$ as $k\rightarrow\infty$. In Section 6, we provide a numerical verification of the theoretical studies. In Section 7, we solve an inverse spectral problem to recover $D$ and the relaxation function from two clusters of eigenvalues. Finally, in the last section, we discuss our results and broader applications.

\section{Formulation of a reduced eigenvalue problem}

In this section, we provide a preliminary study of analyzing the spectrum of operator $A_{N+2}$ that is to introduce a reduced eigenvalue problem with a parameter. This eigenvalue problem can generate a cluster of eigenvalues of $A_{N+2}$.

In the mentioned method to find clustered eigenvalues of $A_{N+2}$, we replace $\xi$ with its generic eigenvalue $-(2k-1)^2$ with $k\in{\mathbb N}$. More precisely, we search for an eigenvalue $\eta$ of $-\partial_x^2$ which gives a non-trivial solution to the following boundary value problem:
\begin{equation}\label{n2.2}
\left\{
\begin{array}{ll}
\partial_x^2u-\eta u=0,\,\, x\in (0,\frac{\pi}{2}),\\
u(0)=0,\,\, \partial_xu(\frac{\pi}{2})=0.
\end{array}
\right.
\end{equation}
It is easy to see that a possible eigenvalue $\eta$ and its associated eigenvector $u$ are given as
\begin{equation}\label{n2.3}
\begin{aligned}
\eta=-(2k-1)^2,\, u=\sin (2k-1)x, \quad k\in \mathbf{N}.
\end{aligned}
\end{equation}

Now, we denote $A_{N+2}^k$ as $A_{N+2}$ with $\xi$ replaced by $-(2k-1)^2$. Then $A_{N+2}^k$ is no longer a differential operator. It is just a multiplication operator on ${\mathbb R}^{N+2}$ by this $(N+2)\times (N+2)$ real matrix $A_{N+2}^k$.
Then, consider the following reduced eigenvalue problem:
\begin{equation}\label{new1.7}
\begin{aligned}
\lambda\,U^k=A_{N+2}^k\,U^k,
\end{aligned}
\end{equation}
where $U^k\in{\mathbb C}^{N+2}$. More precisely, for each $k\in{\mathbb N}$, we look for an eigenvalue $\lambda\not=-r_i,\,1\le i\le N$ of \eqref{new1.7} with the associated eigenfunction $U^k=(u,v,w_1,\cdots,w_n)^{\mathfrak{t}}$ of \eqref{n1.7} given as $\lambda \not= -r_i$.
\begin{equation}\label{n2.4}
\left\{
\begin{array}{ll}
u=\sin (2k-1)x,\\
v=\lambda\sin (2k-1)x,\\
w_i=-(\lambda+r_i)^{-1}(2k-1)^2\sin (2k-1)x, \quad 1\leq i\leq N.
\end{array}
\right.
\end{equation}

It is unclear whether the authors of \cite{LS} were aware of the historical sources of this method because it is a very natural method. The sources of this method that can be identified in \cite{de Verdier}, \cite{Guillemin}, \cite {Gurarie}, \cite{Weinstein}. In these papers, the authors studied the spectrum of the free Laplacian $-\Delta$ perturbed by a real-valued potential $V\in C^\infty(S^d)$ on the $d$-dimensional sphere $S^d$ which is a compact manifold without boundary. It is well known that $-\Delta$ has discrete eigenvalues $\{\mu_\ell=\ell(\ell+d-1)\}_{\ell\in{\mathbb N}\cup\{0\}}$. By perturbing $-\Delta$ by $V$, this eigenvalues are shifted so that $-\Delta+V$ has eigenvalues of the form $$\mu_{\ell m}=\mu_\ell+\zeta_{\ell m},\,\,\ell\in{\mathbb N}\cup\{0\},\,\,1\le m\le p_\ell.$$
For each $\ell\in{\mathbb N}\cup\{0\}$, $\{\mu_{\ell m}\}_{1\le m\le p_\ell}$ is called a cluster of eigenvalues, and clusters of eigenvalues refer to the union of these cluster eigenvalues. The authors of these papers intended to solve the inverse spectral problem identifying $V$ by knowing the clusters of eigenvalues.

Now, note that from \eqref{n1.7}, we have 
\begin{equation}\label{n1.8}
\begin{aligned}
\lambda^2u+D(2k-1)^2 u +  \sum_{i=1}^N b_iw_i=0,
\end{aligned}
\end{equation}
and
\begin{equation}\label{n1.9}
\begin{aligned}
\lambda w_i=-(2k-1)^2 u -r_iw_i,\,\,1\le i\le N,
\end{aligned}
\end{equation}
which are the same to
\begin{equation}\label{n1.10}
\begin{aligned}
(\lambda+r_i)w_i=-(2k-1)^2 u,\,\,\,1\le i\le N.
\end{aligned}
\end{equation}

Multiplying $\Pi_{1\leq j\leq N}(\lambda+r_j)$ to \eqref{n1.8} and using \eqref{n1.10}, we obtain
\begin{equation}\label{n3.1}
\begin{aligned}
(\lambda^2+D(2k-1)^2)\Pi_{1\leq j\leq N}(\lambda+r_j) - (2k-1)^2 \sum_{i=1}^N b_i\Pi_{1\leq j\leq N,j\neq i}(\lambda+r_j)=0,
\end{aligned}
\end{equation}
which is nothing but
\begin{equation}\label{n3.2}
\begin{aligned}
(D+\frac{\lambda^2}{(2k-1)^2})\Pi_{1\leq j\leq N}(\lambda+r_j) -  \sum_{i=1}^N b_i\Pi_{1\leq j\leq N,j\neq i}(\lambda+r_j) =0.
\end{aligned}
\end{equation}
From the arguments deriving \eqref{n3.1}, it is clear that for each $k\in{\mathbb N}$, the roots of \eqref{n3.1} are the eigenvalues of $A_{N+2}$ and the associated eigenvectors are given as
\eqref{n2.4}.

Before closing this section, we define the characteristic polynomial $P^k_N(\lambda)$ associated to $\partial_t-A_{N+2}^k$ and its limit polynomial $P_N(\lambda)$ as $k\rightarrow\infty$ by
\begin{equation}\label{char_poly}
P^k_N(\lambda):=(D+\frac{\lambda^2}{(2k-1)^2})\Pi_{1\leq j\leq N}(\lambda+r_j) -  \sum_{i=1}^N b_i\Pi_{1\leq j\leq N,j\neq i}(\lambda+r_j)
\end{equation}
and
\begin{equation}\label{n3.3}
\begin{aligned}
P_N(\lambda)=D\Pi_{1\leq j\leq N}(\lambda+r_j) -  \sum_{i=1}^N b_i\Pi_{1\leq j\leq N,j\neq i}(\lambda+r_j),
\end{aligned}
\end{equation}
respectively.
In the same way as deriving $P_N^k(\lambda)$, $P_N(\lambda)$ is the charactersite polynomial associated to  
\begin{equation}\label{2.12}
\left\{
\begin{array}{ll}
-D\xi u +  \sum_{i=1}^N b_iw_i=0,\\
w_i'=\xi u -r_iw_i,\,\,1\le i\le N.
\end{array}
\right.
\end{equation}
By assuming each $w_i(0,x)=0$, this is equivalent to the quasi-static equation given as
\begin{equation}\label{2.13}
-D\partial_x^2 u(t,x)+\displaystyle\sum_{i=1}^N b_i\int_0^t\,e^{-r_i(t-\tau)}\partial_x^2 u(\tau,x)\,d\tau=0.
\end{equation}

\section{Roots of characteristic equation and its limit equation}
In this section, we show the existence of all the roots of the characteristic polynomial and the limit polynomials. We start by showing that for the limit polynomial as follows.

\begin{lemma}\label{lem4.1}
There exist exactly $N$ simple real roots $a_j,\,1\le j\le N$ of $P_N(\lambda)$ such that 
$$-r_N<a_N<-r_{N-1}<a_{N-1}<\cdots<-r_2<a_2<-r_1<a_1$$
and
\begin{equation}\label{expression_D1}  
\sum_{i=1}^N \frac{b_i}{a_j+r_i} =D,\,\,1\le j\le N.
\end{equation}
Moreover, we have 
\begin{equation}\label{rep_D}
D=\frac{\sum_{j=1}^N b_j}{\sum_{j=1}^N (a_j+r_j)}.   
\end{equation}
\end{lemma}
\pf. 
For each $1\leq i \leq N$, we have
\begin{equation}\label{n4.1}
\begin{aligned}
P_N(-r_i)&= -   b_i\Pi_{1\leq j\leq N,j\neq i}(-r_i+r_j) \\
&=(-1)^i\cdot |b_i\Pi_{1\leq j\leq N,j\neq i}(-r_i+r_j)|.
\end{aligned}
\end{equation}
Define $f_P(\lambda)$ in $(-r_1,\infty)$ by
\begin{equation}\label{n4.2}
\begin{aligned}
f_P(\lambda)=D-\sum_{i=1}^N b_i(\lambda+r_i)^{-1}.
\end{aligned}
\end{equation}
Since $f_P(\lambda)$ is strictly increasing on $(-r_1,\infty)$, $\lim_{\lambda\rightarrow-r_1+0}f_P(\lambda)=-\infty$, $\lim_{\lambda\rightarrow\infty} f_P(\lambda)\linebreak=D$,
there exists a unique $a_1\in(-r_1,\infty)$ such that 
\begin{equation}\label{n4.3}
\begin{aligned}
f_P(a_1)=0.
\end{aligned}
\end{equation}
Further, for $\lambda\in(-r_1,\infty)$, note that we have 
\begin{equation}\label{n4.4}
\begin{aligned}
P_N(\lambda)=\Pi_{1\leq j\leq N}(\lambda+r_j)f_P(\lambda).
\end{aligned}
\end{equation}
Thus, by \eqref{n4.1}, the intermediate value theorem and the number of roots of $P_N(\lambda)$ is $N$, the roots of $P_N(\lambda)=0$ are
\begin{equation}\label{n4.5}
\begin{aligned}
 a_1, a_2, \cdots,a_{N},
\end{aligned}
\end{equation}
with $r_{j-1}<-a_j<r_{j}$ for $2\leq j \leq N$, and these are the only roots. Then, it is clear that these are all simple roots.
Moreover, comparing the coefficients of $\lambda^{N-1}$ of $P_N(\lambda)=D\,\Pi_{j=1}^N (\lambda-a_j)$ and \eqref{n4.4}, we have 
\begin{equation}\label{n4.6}
\begin{aligned}
D\sum_{j=1}^{N}a_j=\sum_{j=1}^N (b_j-Dr_j),
\end{aligned}
\end{equation}
which immediately yields \eqref{rep_D}.

Finally, to see \eqref{expression_D1} holds for $2\le j\le N$, we only need to note that \eqref{n4.4} holds for $\lambda\in(-r_{j+1},-r_j),\,\,1\le j\le N-1$.
\eproof

In the next lemma, we show the existence of $N$ real roots and $2$ complex conjugate roots of the characteristic polynomial.
\begin{lemma}\label{lem4.2}
There exists at least $N$ real roots of $P_N^k(\lambda)$ such that 
\begin{equation}\label{n4.7}
\begin{aligned}
-r_N<a^k_N<-r_{N-1}<a^k_{N-1}<\cdots<-r_2<a^k_2<-r_1<a^k_1
\end{aligned}
\end{equation}
and
\begin{equation}
\label{expression_D2}
\sum_{i=1}^N \frac{b_i}{a^k_j+r_i} =D+\frac{(a_j^k)^2}{(2k-1)^2},\,\,1\le j\le N.
\end{equation}

The other two roots are contained in the set $B^D_-\cup \{(-r_N,a_1^k)\} $, where $B^D_-:=\{c+id :\frac{-r_N-a_1^k}{2}<c<\frac{-r_1-a_1^k}{2}\} $. 
\end{lemma}
\pf. 
Note that for each $1\leq i \leq N$, we have
\begin{equation}\label{n4.8}
\begin{aligned}
P^k_N(-r_i)&= -   b_i\Pi_{1\leq j\leq N,j\neq i}(-r_i+r_j) \\
&=(-1)^i |b_i\Pi_{1\leq j\leq N,j\neq i}(-r_i+r_j)|.
\end{aligned}
\end{equation}
Define $f^k_P(\lambda)$ in $(-r_1,\infty)$ by
\begin{equation}\label{n4.9}
\begin{aligned}
f^k_P(\lambda)=D+\frac{\lambda^2}{(2k-1)^2}-\sum_{i=1}^N b_i(\lambda+r_i)^{-1}.
\end{aligned}
\end{equation}
Then, for  $\lambda\in(-r_1,\infty)$, we have the expression 
\begin{equation}\label{n4.10}
\begin{aligned}
P^k_N(\lambda)=\Pi_{1\leq j\leq N}(\lambda+r_j)f^k_P(\lambda).
\end{aligned}
\end{equation}
From \eqref{n4.10} and $\Pi_{1\leq j\leq N}(\lambda+r_j)>0$  for  $\lambda\in(-r_1,\infty)$, the roots of $P^k_N(\lambda)=0$ and $f^k_P(\lambda)=0$ in $\lambda$ are the same in $(-r_1,\infty)$. The proof of \eqref{expression_D2} can be done in the same way as in the proof of in Lemma \ref{lem4.1}

Now, observe that
\begin{equation}\label{n4.11}
\begin{aligned}
\lim_{\lambda\rightarrow -r_1+0}f^k_{P}(\lambda)=-\infty.
\end{aligned}
\end{equation}
Further, $f^k_{P}(\lambda)$ is strictly increasing in $(-r_1,\infty)$. Hence, there exists  a unique $a^k_1\in(-r_1,\infty)$ such that 
\begin{equation}\label{n4.12}
\begin{aligned}
f^k_P(a^k_1)=0.
\end{aligned}
\end{equation}
Combining \eqref{n4.8},  \eqref{n4.10}, \eqref{n4.12} and arguing likewise the proof of Lemma \ref{lem4.1}, 
we have $N$ real roots $a_j^k$'s of $P_N^k(\lambda)$ for each $k\in{\mathbb N}$ such that
\begin{equation*}
\begin{aligned}
-r_N<a^k_N<-r_{N-1}<a^k_{N-1}<\cdots<-r_2<a^k_2<-r_1<a^k_1.
\end{aligned}
\end{equation*}
As for \eqref{expression_D2}, it  can be proved in the same way as \eqref{expression_D2} in Lemma \ref{lem4.1}.

Now, for $\lambda<-r_N$, we divide $N$ into two cases. First, we consider the case $N$ is odd. Then $P_N(\lambda)<0$ for  $\lambda<-r_N$.
Hence, for $\lambda<-r_N$,
\begin{equation*}
\begin{aligned}
P^k_N(\lambda)=P_N(\lambda)+\frac{\lambda^2}{(2k-1)^2}\Pi_{1\leq j\leq N}(\lambda+r_j)<0.
\end{aligned}
\end{equation*}

Next, we consider the case $N$ is even. Then $P_N(\lambda)>0$ for  $\lambda<-r_N$. Hence, for $\lambda<-r_N$,
\begin{equation*}
\begin{aligned}
P^k_N(\lambda)=P_N(\lambda)+\frac{\lambda^2}{(2k-1)^2}\Pi_{1\leq j\leq N}(\lambda+r_j)>0.
\end{aligned}
\end{equation*}

Since the degree of $P_N^k(\lambda)$ is $N+2$, $P_N^k(\lambda)$ has two more roots other than $a_j^k,\,1\le j\le N$.
Further, since $P_N^k(\lambda)$ is a polynomial with real coefficients, these two roots are real roots or complex conjugate roots. If $P^k_N(\lambda)$ has other $2$ real roots, then the roots are in the set $(-r_N,a_1^k)$ by the argument given before to determine $a_1^k$.
If the roots are complex roots mentioned above, let the roots of $P^k_N(\lambda)$ be 
\begin{equation*}
\begin{aligned}
a^k_1, a^k_2, a^k_3, \cdots,a^k_{N}, p^k+iq^k,  p^k-iq^k,
\end{aligned}
\end{equation*}
where $r_{j-1}<-a^k_j<r_{j}$  for $2\leq j \leq N$.

Since
\begin{equation}\label{coeff_roots_k}
\begin{aligned}
\sum_{j=1}^{N}a^k_j+2p^k=\sum_{j=1}^N (-r_j),
\end{aligned}
\end{equation}
we obtain that
\begin{equation*}
\begin{aligned}
&-r_N-a_1^k<-r_N-a_1^k+\sum_{j=2}^N (-r_{j-1}-a^k_j)\\
=&2p^k=-r_1-a_1^k+\sum_{j=2}^N (-r_j-a^k_j)<-r_1-a_1^k.
\end{aligned}
\end{equation*}
Thus, we have completed the proof.
\eproof

\section{Asymptotic stability of characteristic roots}\label{sec_asy}

In this section, we first show the locations of the characteristic roots for large $k$ and then their asymptotic behaviors as $k\rightarrow\infty$.
We start by showing the locations of the real characteristic roots in relation to $r_j,\,1\le j\le N$ in the next lemma.

\begin{lemma}\label{lem5.1}
There exists a large positive $k_0$ such that if $k\geq k_0$, then
$P^k_N(\lambda)$ has exactly $N$ simple real roots $a_j^k$, $1\leq j\leq N$ in the interval $(-r_N,a_1^k]$  such that
\begin{equation*}\label{location_realroots}
-r_N<a^k_N<-r_{N-1}<a^k_{N-1}<\cdots<-r_2<a^k_2<-r_1<a^k_1
\end{equation*} 
and other two roots are simple complex conjugate roots in $B^D_-$.
%$B^D_-\cup \{(-r_N,a_1^k]\}$. 
Furthermore 
\begin{equation}\label{limit_real_roots}
\lim_{k\rightarrow\infty}a_j^k=a_j,\,\,1\le j\le N.
\end{equation} 
\end{lemma}
\pf. We write $P_N^k(\lambda)$ in the form:
\begin{equation}\label{n5.1}
\begin{aligned}
P^k_N(\lambda)=&P_N(\lambda)+\frac{\lambda^2}{(2k-1)^2}\Pi_{1\leq j\leq N}(\lambda+r_j)\\
=&P_N(\lambda)+\frac{1}{(2k-1)^2}f(\lambda),
\end{aligned}
\end{equation}
where 
\begin{equation}\label{f_lambda}
f(\lambda):=\lambda^2\cdot\Pi_{1\leq j\leq N}(\lambda+r_j).
\end{equation}

We first show that $P^k_N(\lambda)$ has exactly $N$ real roots in the interval $(-r_N,a_1^k]$.
For that, by Lemma \ref{lem4.1}, we can take small $\delta>0,\,\epsilon>0$ such that
\begin{equation}\label{choice_delta}
\frac{1}{2}\,\min_{2\le j\le N}(a_{j-1}-a_j)>\delta>0,
\end{equation}
\begin{equation}\label{n5.2}
\begin{aligned}
|P'_N(\lambda)|\geq 2\epsilon>0, \quad \lambda \in \cup_{j=1}^N [a_j-\delta,a_j+\delta],
\end{aligned}
\end{equation}
where we have abused the notation "${}'$" to denote the derivative with respect to $\lambda$.
Further, let $m>0$ be a constant such that
\begin{equation}\label{n5.3}
\begin{aligned}
2m\le\min_{\lambda \in [-r_N,a_1+1]\setminus \cup_{j=1}^N [a_j-\delta,a_j+\delta]}|P_N(\lambda)|,
\end{aligned}
\end{equation}
and choose $k_0$ such that 
\begin{equation}\label{n5.4}
\begin{aligned}
|\frac{1}{(2k-1)^2}f(\lambda)|\leq m, \quad |\frac{1}{(2k-1)^2}f'(\lambda)|\leq \epsilon, \quad \lambda \in  [-r_N,a_1+1].
\end{aligned}
\end{equation}
hold for $k\ge k_0$.
It is easy to see that
\begin{equation}\label{n5.5}
\begin{aligned}
|P^k_N(\lambda)|=&|P_N(\lambda)+\frac{1}{(2k-1)^2}f(\lambda)|\\
\geq&|P_N(\lambda)|-|\frac{1}{(2k-1)^2}f(\lambda)|\geq m, \quad \lambda \in [-r_N,a_1+1]\setminus \cup_{j=1}^N [a_j-\delta,a_j+\delta]
\end{aligned}
\end{equation}
and
\begin{equation}\label{n5.6}
\begin{aligned}
|(P^k_N)'(\lambda)|=&|P'_N(\lambda)+\frac{1}{(2k-1)^2}f'(\lambda)|\\
\geq&|P'_N(\lambda)|-|\frac{1}{(2k-1)^2}f'(\lambda)|\geq \epsilon, \quad \lambda \in \cup_{j=1}^N [a_j-\delta,a_j+\delta].
\end{aligned}
\end{equation}
From Lemma \ref{lem4.1} and \eqref{n5.2}, we have
\begin{equation}\label{n5.7}
\begin{aligned}
P_N(a_j-\delta)\cdot P_N(a_j+\delta)<0, \quad  1\leq j \leq N.
\end{aligned}
\end{equation}
Combining \eqref{n5.3}, \eqref{n5.5} and \eqref{n5.7}, we obtain that
\begin{equation}\label{n5.8}
P^k_N(a_j-\delta)\cdot P^k_N(a_j+\delta)<0, \quad  1\leq j \leq N.
\end{equation}
From \eqref{n5.6} and \eqref{n5.8}, there exists a single root $a^k_j$ in
$(a_j-\delta,a_j+\delta)$ for each $1\leq j \leq N$, which yields \eqref{limit_real_roots}. As for the statement of two other roots in addition to $a_j^k,\,1\le j\le N$, it has already been proved in Lemma \ref{lem4.2}.
Thus, we have completed the proof.
\eproof

\begin{comment}
If $k$ is large enough, then
the roots of $P^k_N(\lambda)$ will be
\begin{equation}\label{n5.9}
\begin{aligned}
a^k_1, a^k_2, a^k_3, \cdots,a^k_{N}, p^k+iq^k,  p^k-iq^k,
\end{aligned}
\end{equation}
where $r_{j-1}<-a^k_j<r_{j}$ and $\lim_{k\rightarrow \infty}a^k_j=a_j$ for $1\leq j \leq N$.
\end{comment}

\begin{comment}From \eqref{n3.2}, we have
\begin{equation}\label{n5.10}
\begin{aligned}
\sum_{j=1}^{N}a^k_j+2p^k=\sum_{j=1}^N (-r_j).
\end{aligned}
\end{equation}
Also, from \eqref{n3.3}, we have
\begin{equation}\label{n5.11}
\begin{aligned}
\sum_{j=1}^{N}a_j=\sum_{j=1}^N (-r_j)+\frac{1}{D}(\sum_{j=1}^N b_j).
\end{aligned}
\end{equation}
\end{comment}

In the following two lemmas, we give estimates of the convergence of sequences $\{a_j^k\}$ and $\{(p^k,\,q^k)\}$. As for $\{a_j^k\}$, it is given as follows.
\begin{lemma}\label{lem5.2}
Let $k\geq k_0$ and $a_j^k$, $1\leq j\leq N$ be $N$
   real roots of $P^k_N(\lambda)$. Then, we have
\begin{equation}\label{n5.12}
\left\{
\begin{array}{ll}
-r_j<a^k_j<a_j\\
0\leq a_j-a_j^k\leq \frac{M_1}{k^2},
\end{array}
\right.
\end{equation} 
where $M_1$ is a positive constant independent of $k\in{\mathbb N}$.
\end{lemma}
\pf. 
We rewrite $P_N^k(\lambda)$ as follows:
\begin{equation}\label{n5.13}
\begin{aligned}
P^k_N(\lambda)=&P_N(\lambda)+\frac{\lambda^2}{(2k-1)^2}\Pi_{1\leq j\leq N}(\lambda+r_j)\\
=&D\cdot\Pi_{1\leq j\leq N}(\lambda-a_j)+\frac{1}{(2k-1)^2}f(\lambda)
\end{aligned}
\end{equation}
with $f(\lambda)$ given by \eqref{f_lambda}. Here, from Lemma \ref{lem5.1}, there exists a single root $a^k_j$ in
$(a_j-\delta,a_j+\delta)$ for each $1\leq j \leq N$. Also, by $a_j>a_{j+1},\,1\le j\le N-1$ and \eqref{choice_delta}, we have
\begin{equation}\label{n5.14}
\begin{aligned}
0=&P^k_N(a_l^k)=D\cdot\Pi_{1\leq j\leq N}(a_l^k-a_j)+\frac{1}{(2k-1)^2}f(a_l^k)\\
=&(-1)^{l+1}(a_l^k-a_l)|D\cdot\Pi_{1\leq j\leq N, j\neq l}(a_l^k-a_j)|+(-1)^{l+1}\frac{1}{(2k-1)^2}|f(a_l^k)|
\end{aligned}
\end{equation}
for each $1\leq l\leq N$. 
This implies that for each $1\leq l\leq N$, we have
\begin{equation}\label{n5.15}
\begin{aligned}
a_l^k-a_l=(-1)|D\cdot\Pi_{1\leq j\leq N, j\neq l}(a_l^k-a_j)|^{-1}\cdot\frac{1}{(2k-1)^2}|f(a_l^k)|.
\end{aligned}
\end{equation}
Therefore, 
\begin{equation}\label{n5.16}
\begin{aligned}
0\leq a_l-a_l^k\leq \frac{M_1}{k^2}
\end{aligned}
\end{equation}
holds with a  positive constant $M_1$ independent of $k\in{\mathbb N}$.
\eproof

\begin{lemma}\label{lem5.3}
Let $k\geq k_0$ and $p^k+iq^k$, $p^k-iq^k$ with $p_k, q_k\in{\mathbb R}$ be $2$
complex roots of $P^k_N(\lambda)$. Then, we have
\begin{equation}\label{n5.17}
\left\{
\begin{array}{ll}
|p^k+\frac{1}{2D}(\sum_{j=1}^N b_j)|\leq \frac{M_2}{k^2},\\
|q^k-(2k-1)\sqrt{D}|\leq  M_3(2k-1)^{-1},
\end{array}
\right.
\end{equation} 
where $M_2$ and $M_3$ are positive constants independent of $k\in{\mathbb N}$.
\end{lemma}
\pf. For the first inequality of \eqref{n5.17}, it follows by combining \eqref{n4.6}, \eqref{coeff_roots_k} and \eqref{n5.12} by taking $M_2=2^{-1}M_1$. We proceed as follows to prove the second inequality of \eqref{n5.17}.
Using the roots of $P^k_N(\lambda)$, we rewrite $(2k-1)^2P^k_N(\lambda)$ to have
\begin{equation}\label{n5.18}
\begin{array}{ll}
(\lambda^2+D(2k-1)^2)\Pi_{1\leq j\leq N}(\lambda+r_j) - (2k-1)^2 \sum_{i=1}^N b_i\Pi_{1\leq j\leq N,j\neq i}(\lambda+r_j)\\
\qquad\qquad\quad=\Pi_{1\leq j\leq N}(\lambda-a_j^k)\big(\lambda^2-2p^k\lambda+(p^k)^2+(q^k)^2\big).
\end{array}
\end{equation}
Comparing the coefficients $\lambda^N$ of \eqref{n5.18}, we obtain that
\begin{equation}\label{n5.19}
\begin{array}{ll}
\Pi_{1\leq j<l\leq N}(r_jr_l) +D (2k-1)^2\\
\qquad\qquad\qquad=(p^k)^2+(q^k)^2+2p^k(\sum_{j=1}^N a_j^k)+\Pi_{1\leq j<l\leq N}(a_j^ka_l^k).
\end{array}
\end{equation}
Combining this with the first inequality of \eqref{n5.17}, we have
\begin{equation}\label{n5.20}
\begin{aligned}
|(q^k)^2-D(2k-1)^2|\leq M_2.
\end{aligned}
\end{equation}
For $k\ge2^{-1}\{1+\sqrt{2M_3(1-\kappa)}\}$ with $M_3=2^{-1} D^{-1/2}M_2,\,0<\kappa\ll 1$,   \eqref{n5.20} implies 
\begin{equation}\label{n5.21}
\begin{aligned}
&\sqrt{D}-M_3(2k-1)^{-2}\leq \sqrt{D-M_2(2k-1)^{-2}}\leq (q^k)(2k-1)^{-1}\\
\leq &\sqrt{D+M_2(2k-1)^{-2}}\leq \sqrt{D}+M_3(2k-1)^{-2},
\end{aligned}
\end{equation}
which is equivalent to
\begin{equation}\label{n5.22}
\begin{aligned}
|(q^k)-(2k-1)\sqrt{D}|\leq  M_3(2k-1)^{-1}.
\end{aligned}
\end{equation}
Thus, we proved the second inequality of \eqref{n5.17}.
\eproof

We closed this section by giving upper estimates of the distances of adjacent real eigenvalues of $\{a_j^k\}$ for each $j$ and complex eigenvalues of $\{p^k+iq^k\}$ which immediately follow from \eqref{n5.12} and \eqref{n5.17}. 

\begin{remark} For each fix $1\le j\le N$, the following estimates hold for $k\ge k_0${\,\rm :}
\begin{equation}\label{dist_ev}
\left\{
\begin{array}{ll}
|a_j^{k+1}-a_j^k|\le 2\frac{M_1}{k^2},\,\, |p^{k+1}-p^k|\le 2\frac{M_2}{k^2},\\
|q^{k+1}-q^k|\le 2\frac{M_3}{2k-1}+2\sqrt{D}.
\end{array}
\right.
\end{equation}
\end{remark}

\section{Explicit form of $k_0$ in Lemma \ref{lem5.1}  }
In this section, we give an explicit form of $k_0$. 
To begin with, we estimate $a_1$ as in the following lemma.
\begin{lemma}\label{lem6.1}
 $a_1\leq \frac{NB}{D}$, where $B=\max_{1\leq i\leq N}b_i$.
\end{lemma}
\pf. By Lemma \ref{lem4.1}, we have $a_1\in (-r_1,\infty)$ and $\sum_{i=1}^N \frac{b_i}{a_1+r_i} =D$. Since the statement is trivial for $a_1\le0$, we only need to consider the case $a_1>0$. Then, by $0<r_1<\cdots<r_N$, we have
\begin{equation}\label{n6.3}
\begin{aligned}
D\leq B\sum_{i=1}^N \frac{1}{a_1+r_i}\leq NB\frac{1}{a_1+r_1}.
\end{aligned}
\end{equation}
Hence, we have
\begin{equation}\label{n6.4}
\begin{aligned}
a_1<a_1+r_1\leq \frac{NB}{D}.
\end{aligned}
\end{equation}
\eproof

\begin{lemma}\label{lem6.2}
Let 
\begin{equation}\label{b_r_mu}
\left\{
\begin{array}{ll}
b:=\min_{1\leq i\leq N}b_i,\,\, r:=\min_{1\leq i\leq N-1}|r_{i+1}-r_i|,\\
\\
\mu:=\min\{\frac{br^{N-1}}{4D(2r_N)^{N-1}+4NB(2r_N)^{N-2}+1},\frac{r}{4},\frac{r_1}{4}\}.
\end{array}
\right.
\end{equation}
Then, we have
\begin{equation}\label{n6.5}
\begin{aligned}
|P_N(\lambda)|\geq \frac{br^{N-1}}{2} >0,\,\,\,\lambda\in \cup_{j=1}^N [-r_j-\mu,-r_j+\mu],
\end{aligned}
\end{equation}
which implies 
\begin{equation}\label{dist_of_roots}
|a_i-a_j|\geq 2\mu,\,\,1\leq i <j \leq N.    
\end{equation}
\end{lemma}
\pf. Fix any $1\le l\le N$ and let $\lambda \in  [-r_l-\mu,-r_l+\mu]$. By using \eqref{n3.3}, we estimate $|P_N(\lambda)|$ from below to get \eqref{n6.5} as follows: 
\begin{equation*}
\begin{aligned}
|P_N(\lambda)|=&|-b_l\Pi_{1\leq j\leq N,j\neq l}(\lambda+r_j)+D\Pi_{1\leq j\leq N}(\lambda+r_j) -  \sum_{1\leq i\leq N,i\neq l} b_i\Pi_{1\leq j\leq N,j\neq i}(\lambda+r_j)|\\
\geq & \frac{3br^{N-1}}{4}-|D\Pi_{1\leq j\leq N}(\lambda+r_j)|-| \sum_{1\leq i\leq N,i\neq l} b_i\Pi_{1\leq j\leq N,j\neq i}(\lambda+r_j)|\\
\geq & \frac{3br^{N-1}}{4}-\mu D(2r_N)^{N-1}-\mu NB(2r_N)^{N-2}\\
\geq & \frac{br^{N-1}}{2}.
\end{aligned}
\end{equation*}
\eproof

\begin{lemma}\label{lem6.3}
Let $R:=2r_N+\frac{NB}{D}+1$, and take $\delta$ in \eqref{choice_delta} as $\delta:=\frac{r^{N-2}}{2NR^{N-2}}\mu$ which is possible by \eqref{dist_of_roots}.
Then the constants $\epsilon$ in \eqref{n5.2} and $m$ in \eqref{n5.3} can be taken $\epsilon=D\mu r^{N-2}/4$ and $m=D\mu\delta r^{N-2}/2$, respectively.
\end{lemma}
\pf. 
Since
\begin{equation}\label{n6.1}
\begin{aligned}
P_N(\lambda)=D\cdot\Pi_{1\leq j\leq N}(\lambda-a_j)
\end{aligned}
\end{equation}
by Lemma \ref{lem4.1},
we have
\begin{equation}\label{n6.2}
\begin{aligned}
P_N'(\lambda)=&\sum_{l=1}^N D\cdot\Pi_{1\leq j\leq N,j\neq l}(\lambda-a_j)\\
=& D\cdot\Pi_{1\leq j\leq N,j\neq i}(\lambda-a_j)+\sum_{1\leq l\leq N,l\neq i} D\cdot\Pi_{1\leq j\leq N,j\neq l}(\lambda-a_j)
\end{aligned}
\end{equation}
for each $1\le i\le N$. From \eqref{n6.2} and \eqref{n6.5}, we have that
\begin{equation}\label{n6.7}
\begin{aligned}
|P_N'(\lambda)|\geq& |D\cdot\Pi_{1\leq j\leq N,j\neq i}(\lambda-a_j)|-\sum_{1\leq l\leq N,l\neq i} |D\cdot\Pi_{1\leq j\leq N,j\neq l}(\lambda-a_j)|\\
\geq& D\mu r^{N-2}-DN\delta R^{N-2}\\
\geq& D\mu r^{N-2}/2
\end{aligned}
\end{equation}
for $\lambda \in  [a_i-\delta,a_i+\delta]$ and $1\leq i \leq N$.

From \eqref{n6.1}, we have for $\lambda \in [-r_N,a_1+1]\setminus \cup_{j=1}^N [a_j-\delta,a_j+\delta]$ that
\begin{equation}\label{n6.8}
\begin{aligned}
|P_N(\lambda)|= |D\cdot\Pi_{1\leq j\leq N}(\lambda-a_j)|
\geq D\mu \delta r^{N-2}.
\end{aligned}
\end{equation}
Hence, we can take the mentioned $m$. 
\eproof
\begin{lemma}\label{lem6.4}
Let $k_0=\frac{R^{(N+2)/2}}{\sqrt{m}}+\frac{R^{(N+1)/2}}{\sqrt{\epsilon}}+1$, where $m$ and $\epsilon$ are defined in Lemma \ref{lem6.3}. Then \eqref{n5.4} holds for $k\geq k_0$. 
\end{lemma}
\pf. By \eqref{f_lambda} and $R$ given in Lemma \ref{lem6.3}, we easily have $|f(\lambda)|\leq R^{N+2}$ and $|f'(\lambda)|\leq (N+1)R^{N+1}$ for $\lambda \in  [-r_N,a_1+1]$. Then, by recalling the proof of Lemma \ref{lem5.1}, we can have \eqref{n5.4} for $k\geq k_0$.
\eproof

\section{Numerical verification of theoretical results}

In this section, we numerically verify the accuracy of our spectral analysis. More precisely, by fixing $h,\,N,\,\text{$r_i$'s}$ to $h=1$, $N=5\,\,\text{or}\,\,9$, $r_i=5*i,\,i=1,\,2,\,\cdots,\,N$, we numerically calculate the roots $a_j^k$ of $P_N^k$ defined by \eqref{char_poly} for different $k, N$ for the cases $D>h,\,D=h,\,D<h$. We also do the same for the roots $a_j$ of $P_N$ defined by \eqref{n3.3}.  Then, we show the numerical convergence of the roots of $P_N^k$ proved in Lemmas \ref{lem5.2}--\ref{lem5.3} as $k$ increases.

\begin{figure}[htp]
\centering
\begin{tabular}{ll}
(a) real roots for $D=0.5$ &(b) complex roots for $D=0.5$  \\
\includegraphics[width=0.5\textwidth]{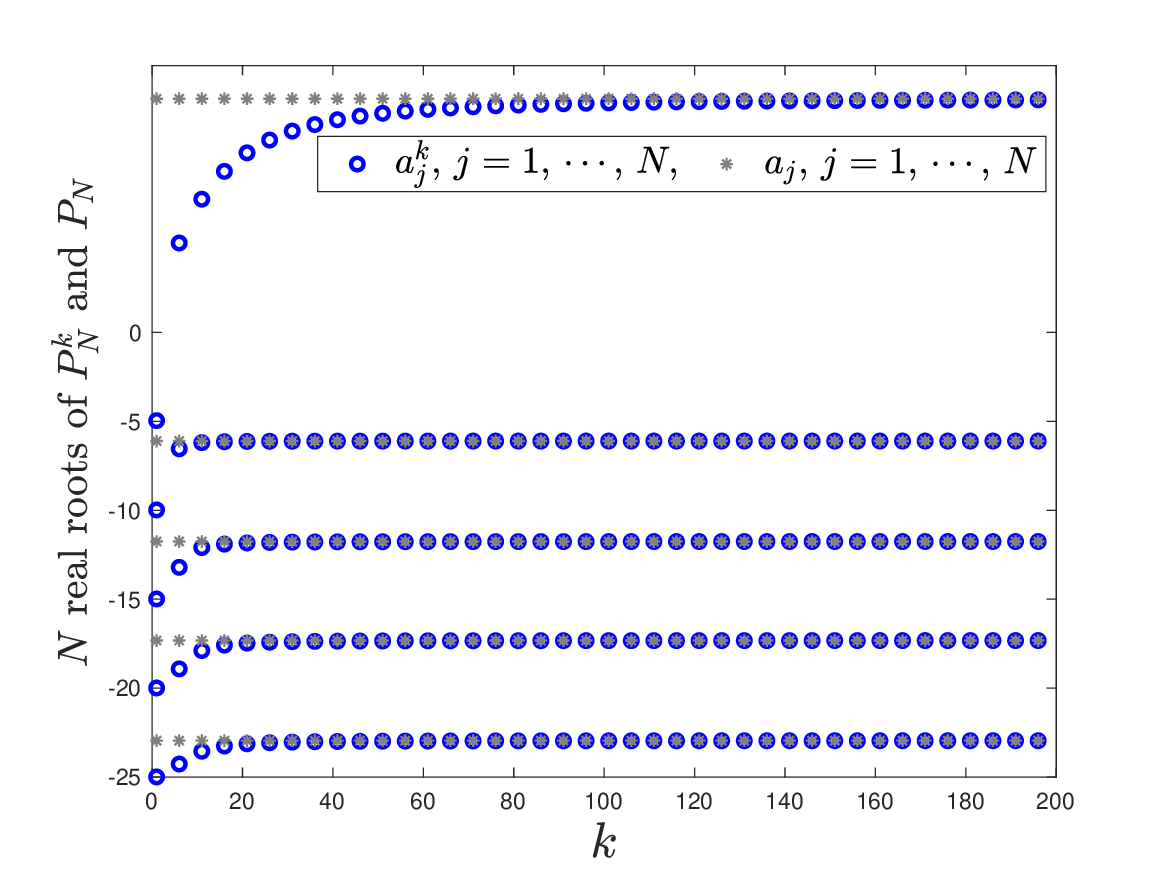}
&\includegraphics[width=0.5\textwidth]{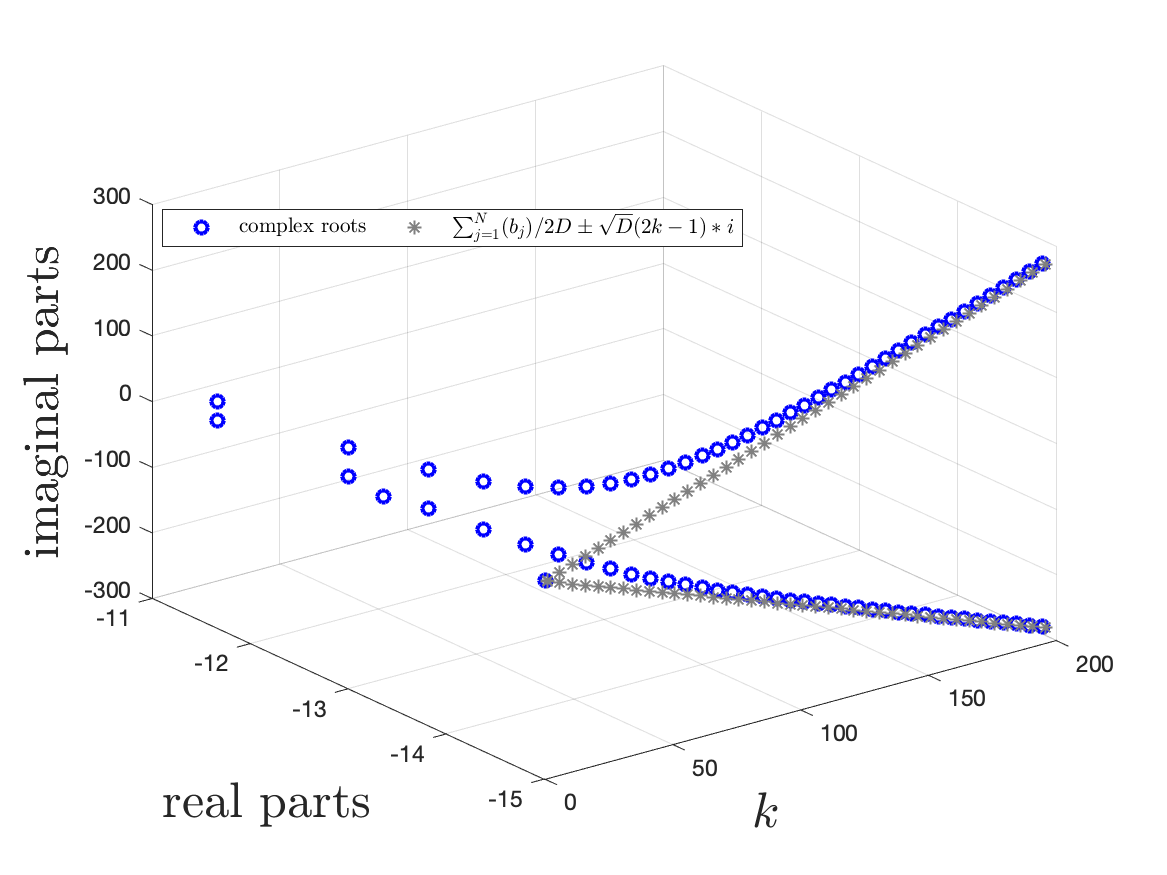}
%&\includegraphics[width=0.33\textwidth]{pic/complex_roots_iN5.eps}
\\
(c) real roots for $D=1$ &(d) complex roots for $D=1$ \\
\includegraphics[width=0.5\textwidth]{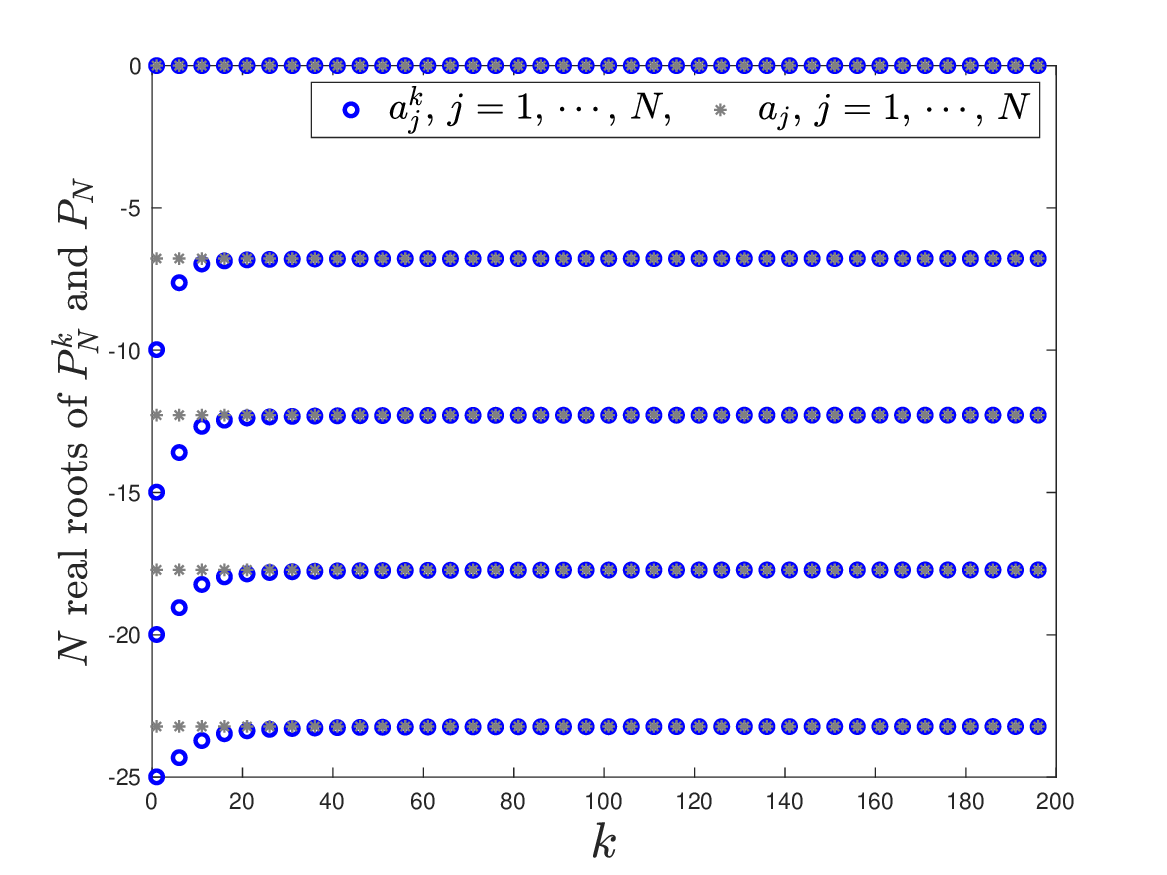}
&\includegraphics[width=0.5\textwidth]{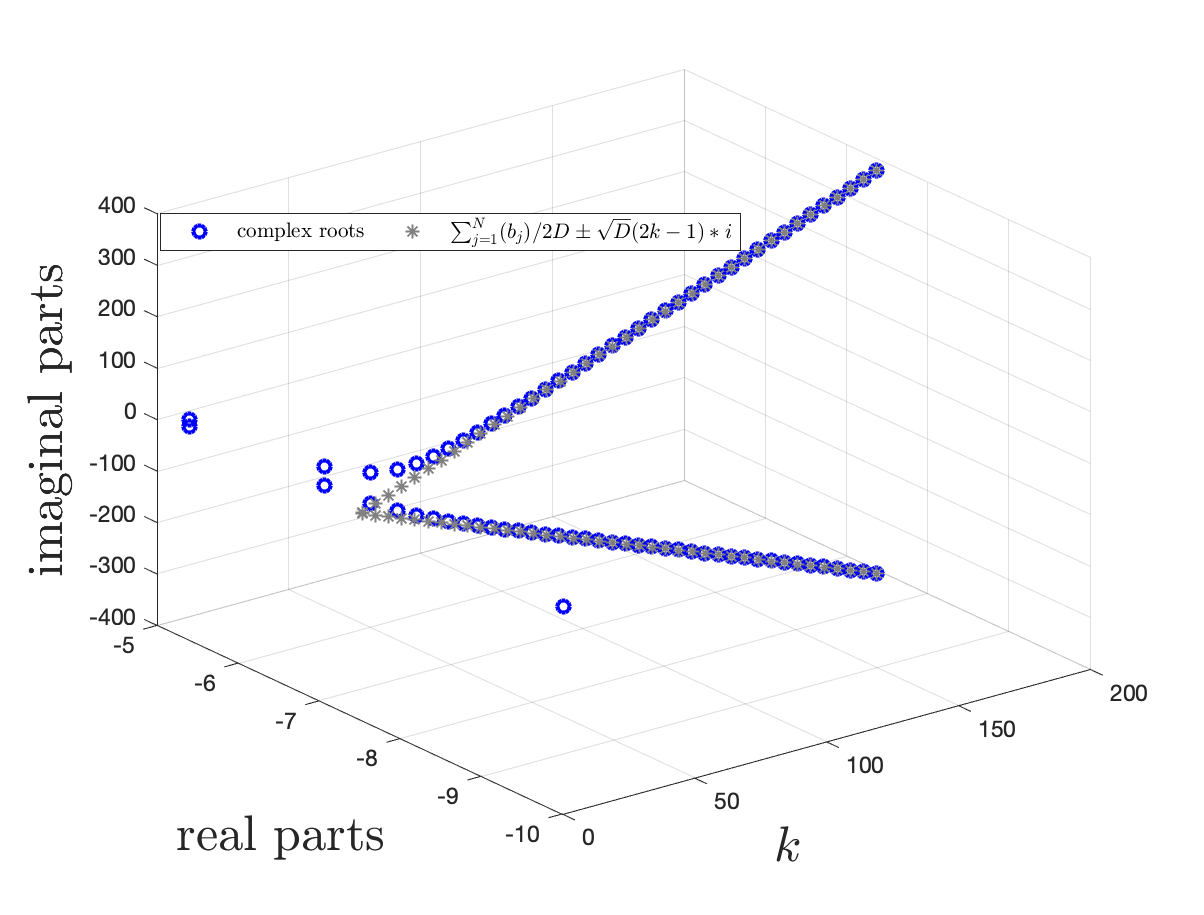}
\\
(e) real roots for $D=5$ &(f) complex roots for $D=5$ \\
\includegraphics[width=0.5\textwidth]{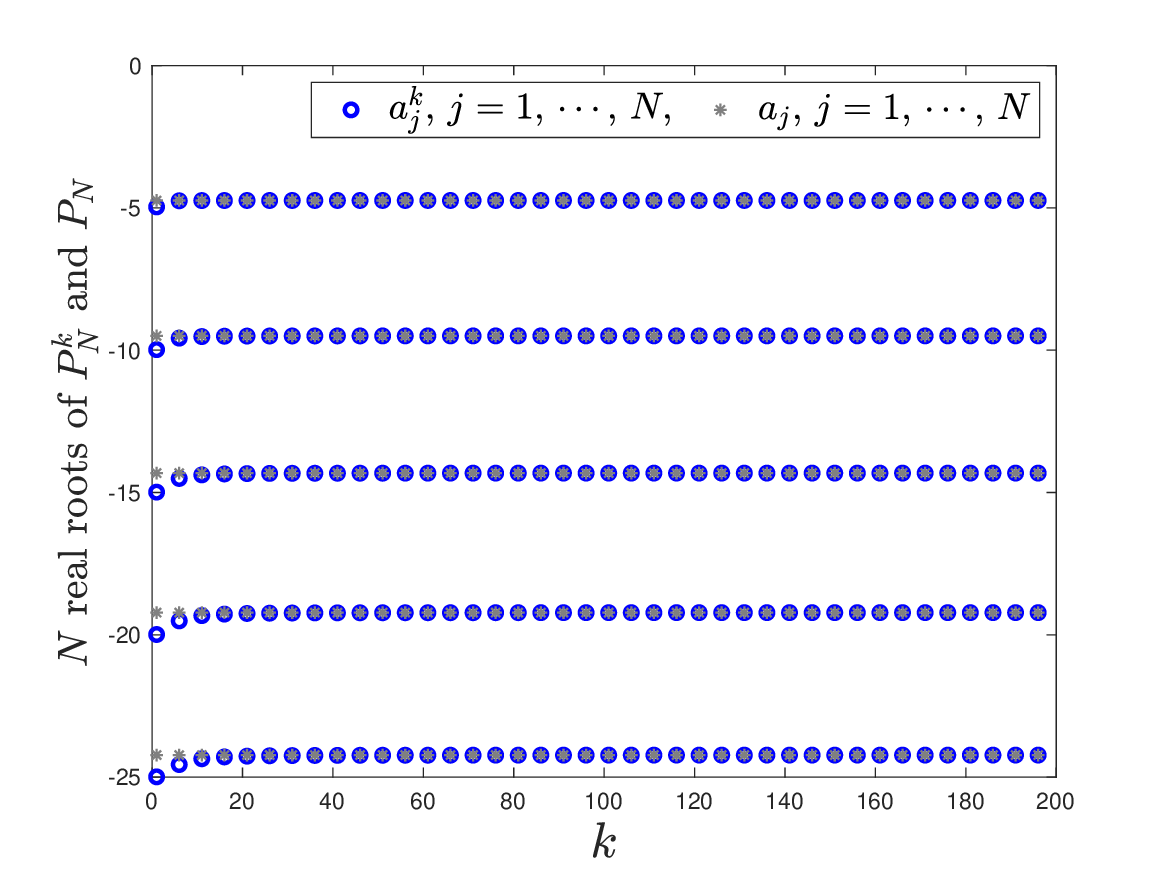}
&\includegraphics[width=0.5\textwidth]{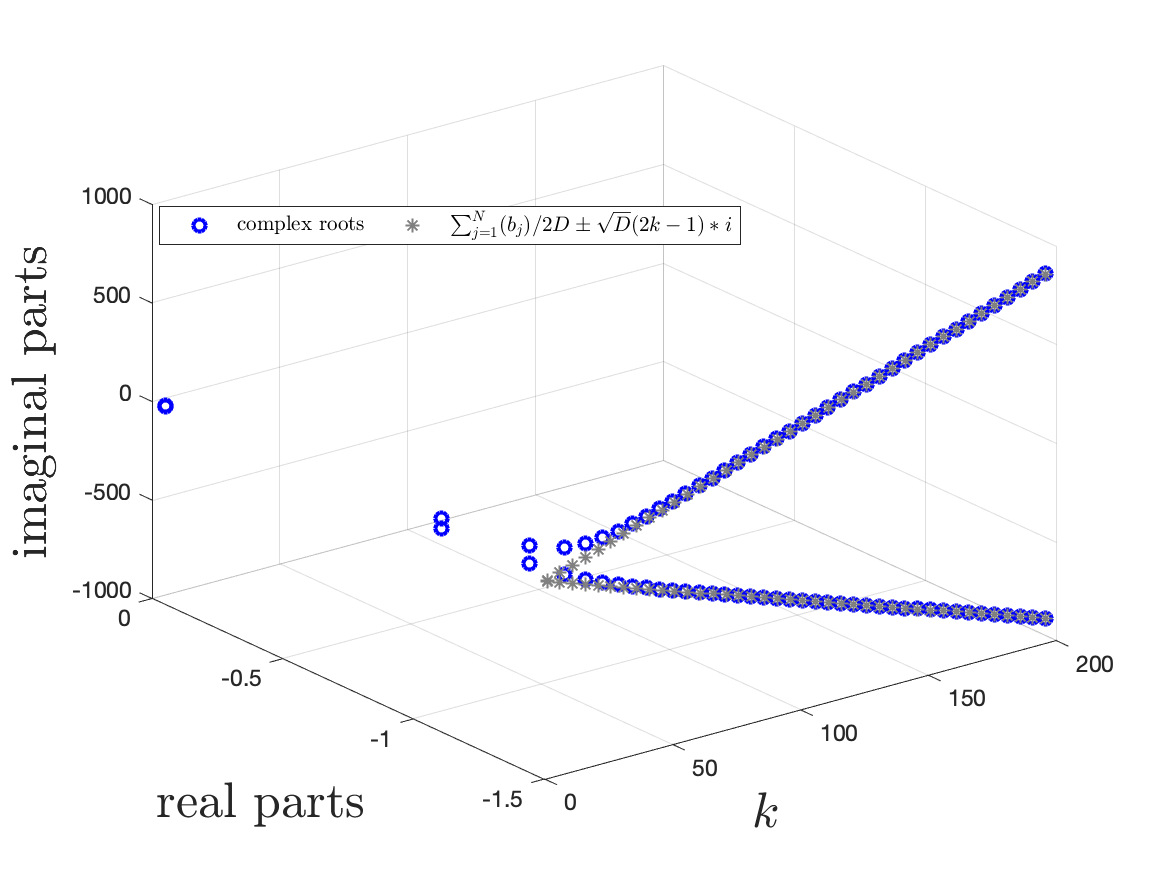}
\end{tabular}
\caption{Roots of $P_N^k$ and $P_N$ for $N=5$ and $D=0.5,\,1,\,5$}
\label{fig_N5g05_1_5}
\end{figure}

\begin{example}\label{exm_1}
Let $N=5$. The roots of $P_N^k$ and $P_k$ are calculated and plotted in Figure \ref{fig_N5g05_1_5}. 
\end{example}

In Figure \ref{fig_N5g05_1_5}, for each $D$ and $k$, there are $5$ real roots for both $P_N^k,\,P_N$, and additional $2$ complex roots for $P_N^k$ complex conjugate to each other. For the real roots, we can observe that the locations of real roots $a_j^k,\,j=1,\,2,\,\cdots,5$ satisfy \eqref{location_realroots}. For the each mentioned case on $D$, the real roots $a_j^k,\,j=1,\,2,\,\cdots,5$ of $P_N^k$ approach the corresponding roots $a_j,\,j=1,\,2,\,\cdots,5$ of $P_N$ as $k$ increases. As for the complex roots $p_k\pm i q_k$ of $P_N^k$, they approach $\sum_{j=1}^N\frac{b_j}{2D}\pm\sqrt{D}(2k-1)$ as $k$ increases.

\begin{figure}[htp]
\centering
\begin{tabular}{ll}
(a) real roots for $D=0.5$ &(b) complex roots for $D=0.5$  \\
\includegraphics[width=0.5\textwidth]{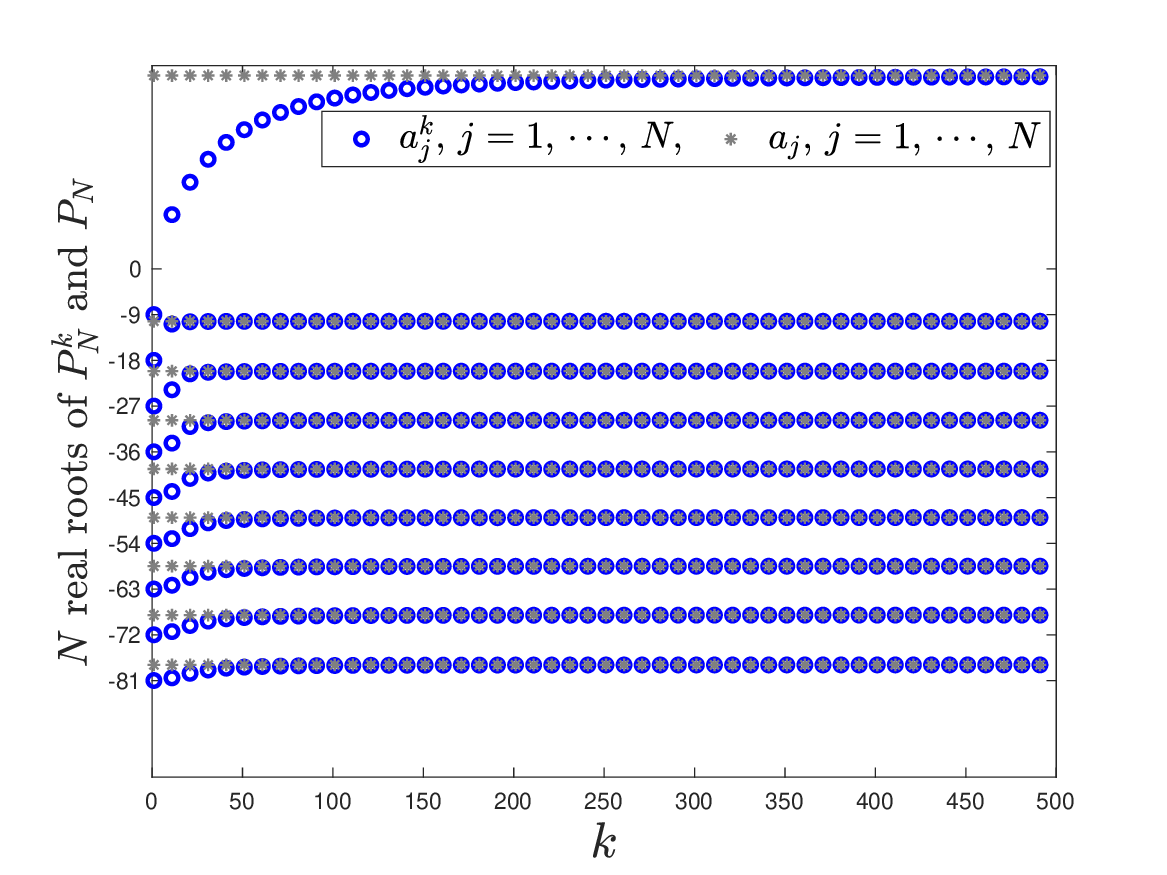}
&\includegraphics[width=0.5\textwidth]{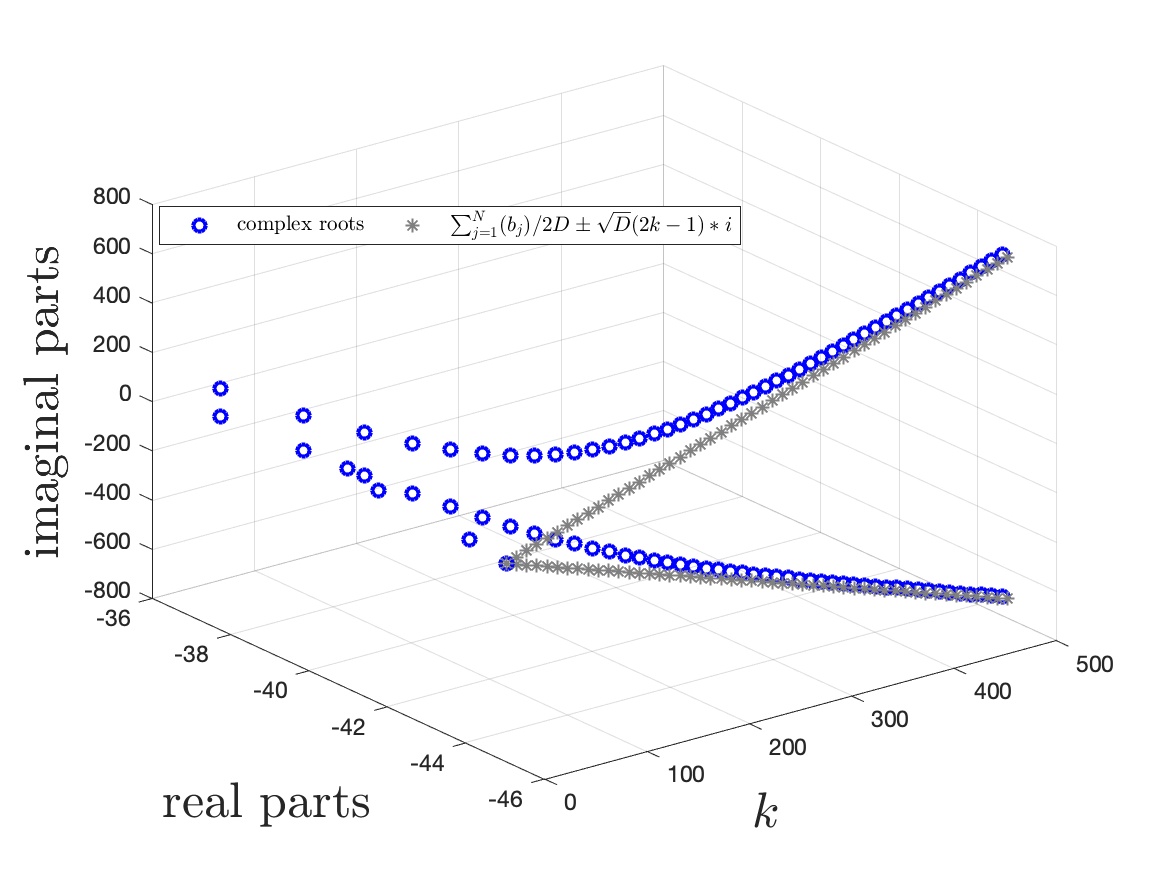}
\\
(c) real roots for $D=1$ &(d) complex roots for $D=1$  \\
\includegraphics[width=0.5\textwidth]{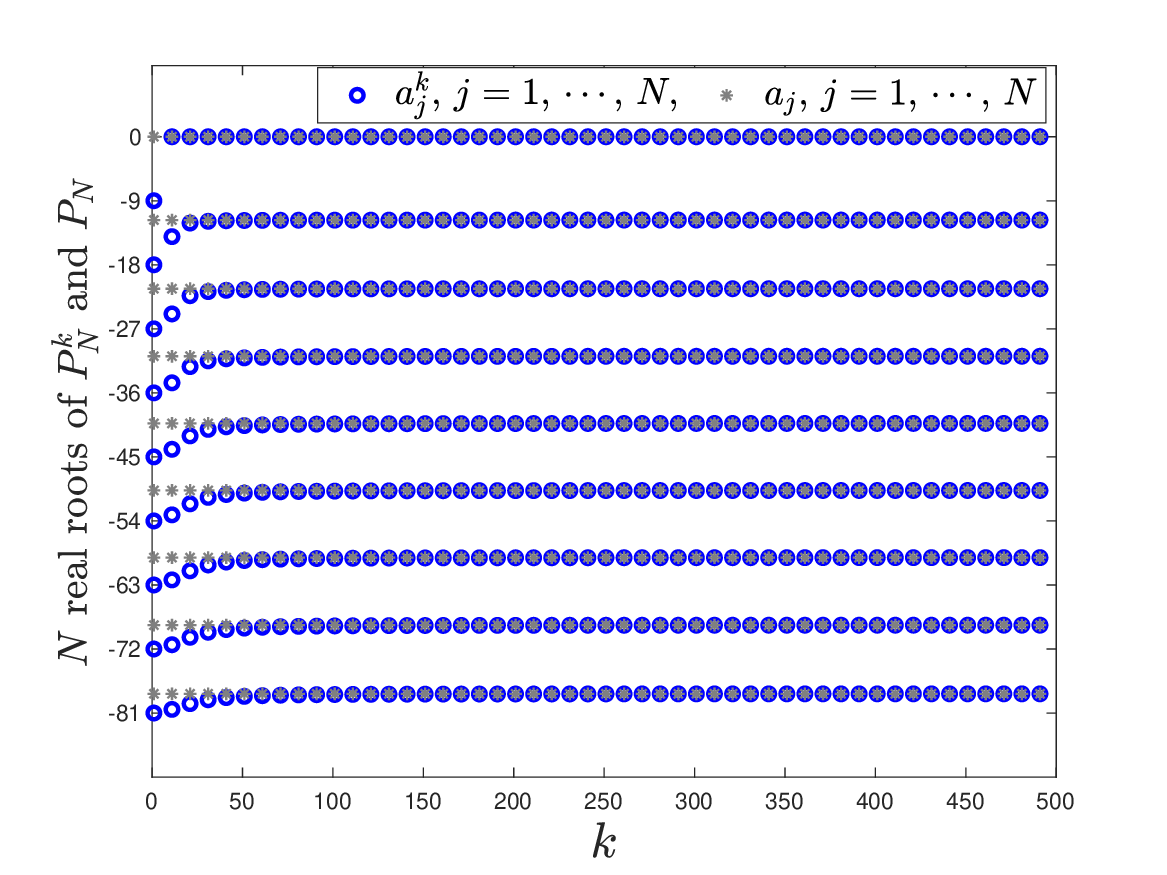}
&\includegraphics[width=0.5\textwidth]{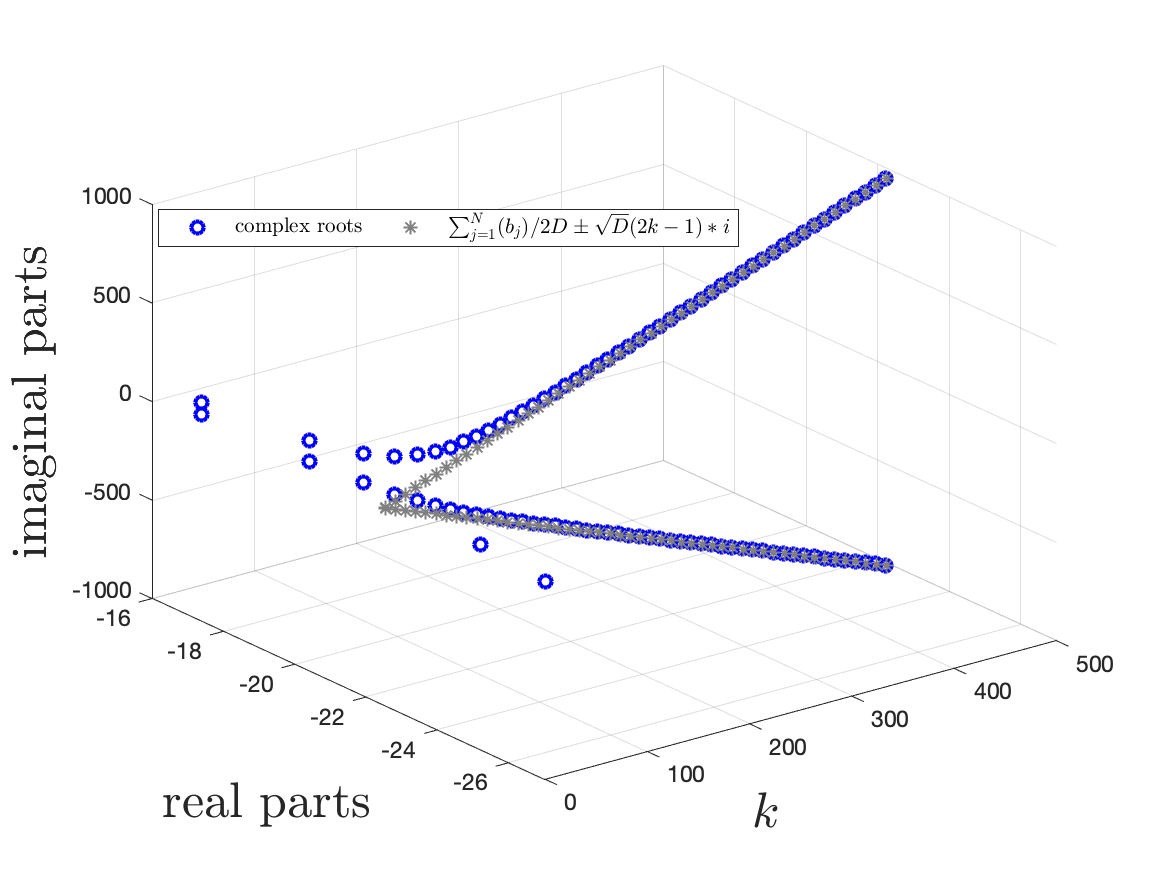}\\
(e)  real roots for $D=5$ &(f) complex roots for $D=5$ \\
\includegraphics[width=0.5\textwidth]{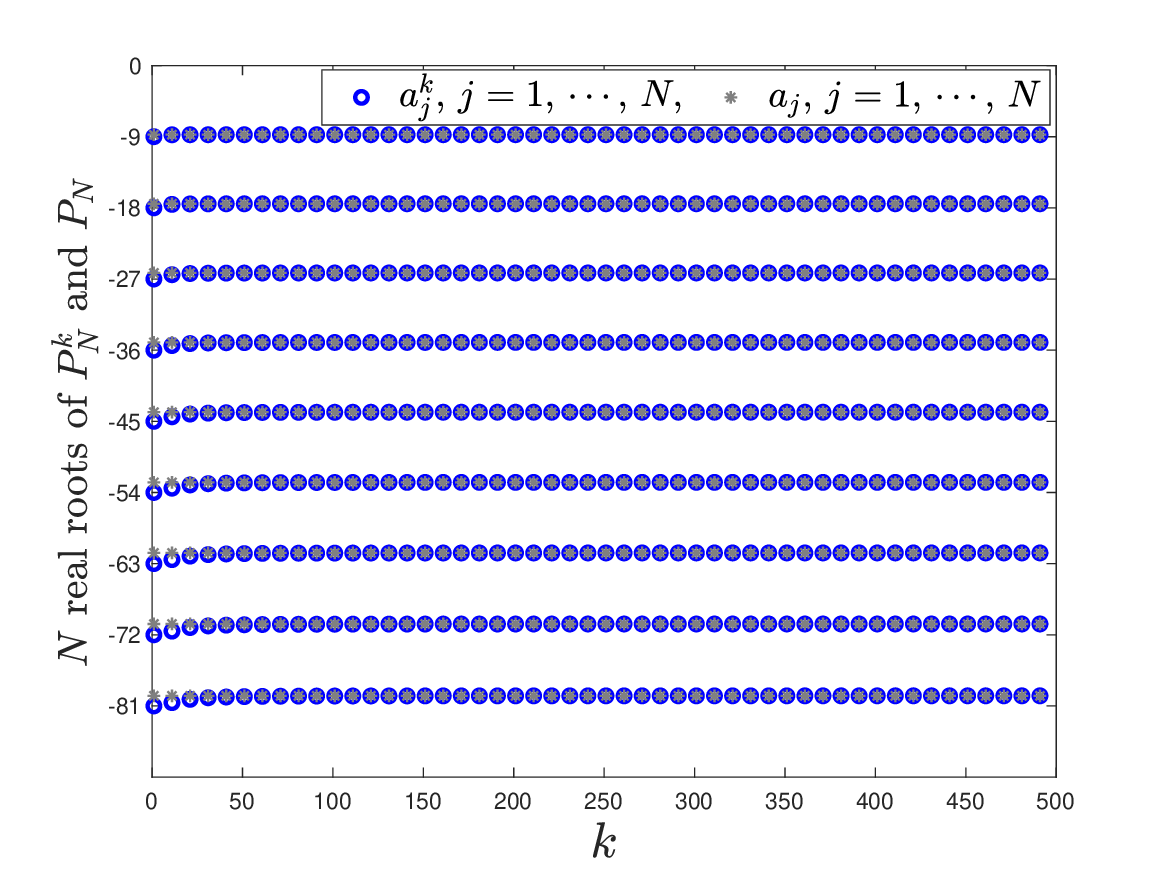}
&\includegraphics[width=0.5\textwidth]{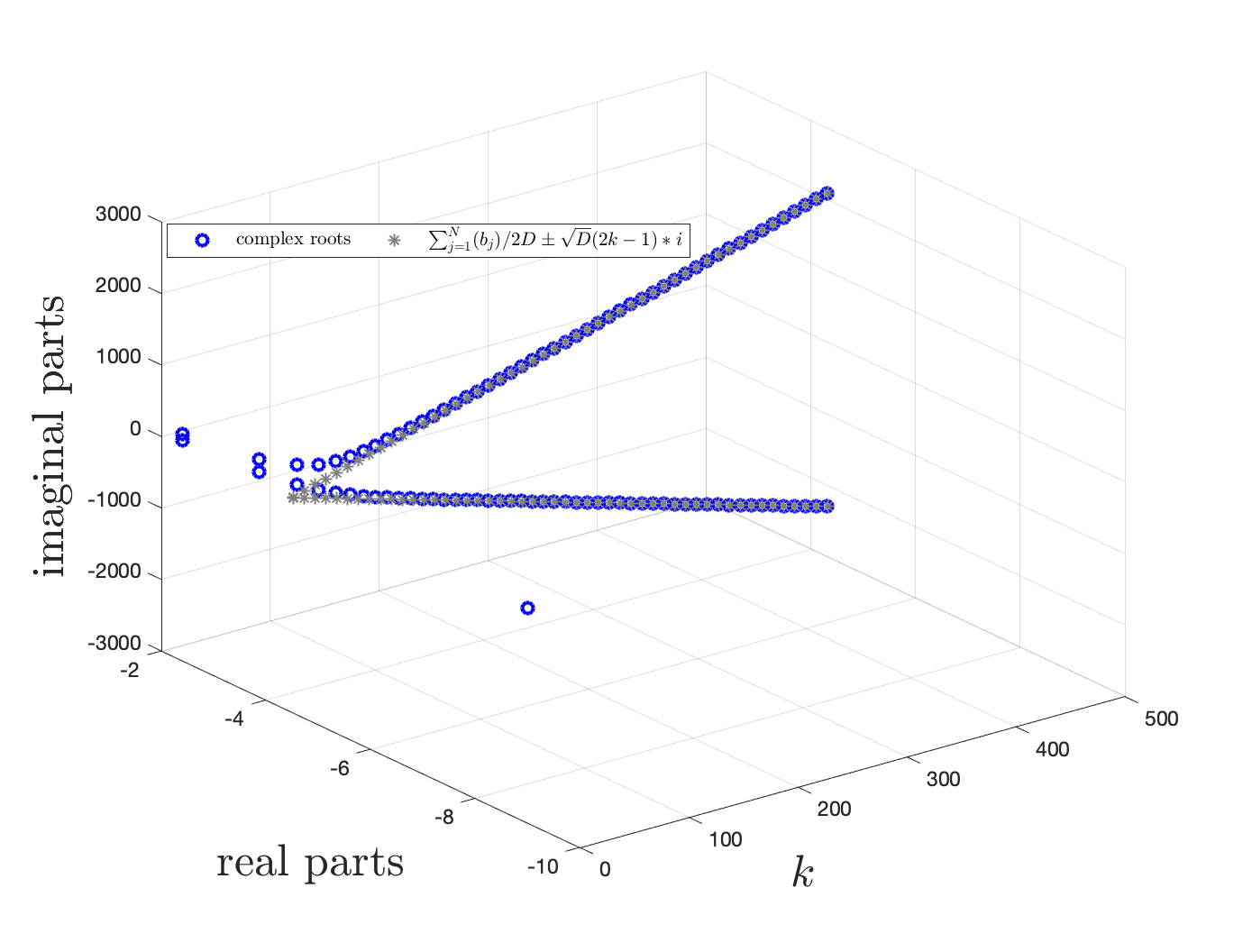}
\end{tabular}
\caption{Roots of $P_N^k$ and $P_N$ for $N=9$ and $D=0.5,\,1,\,5$}
\label{fig_N9g05_1_5}
\end{figure}

\begin{example}\label{exm_2}
Let $N=9$. The roots of $P_N^k$ and $P_k$ are calculated and plotted in Figure \ref{fig_N9g05_1_5}. 
\end{example}

%In the following, we numerically tested the asymptotic behavior of roots of $P_N^k(\lambda)$. The results for $N=5,\,9$ and $D=1,\,5$ are shown in Figures \ref{fig_N5}--\ref{fig_N9_g5}.

In Figure \ref{fig_N9g05_1_5}, we can observe profiles of roots similar to \ref{exm_1} except for their numbers. However, the convergence rates of the roots of $P_N^k$ are slower than those of Example \ref{exm_1}, especially for the cases $D=0.5$ and $D=1$. 

\begin{figure}[htp]
\centering
\begin{tabular}{lll}
(a) real roots for $D=0.5$ &(b) real part for $D=0.5$ &(c) imag. part for $D=0.5$  \\
\includegraphics[width=0.33\textwidth]{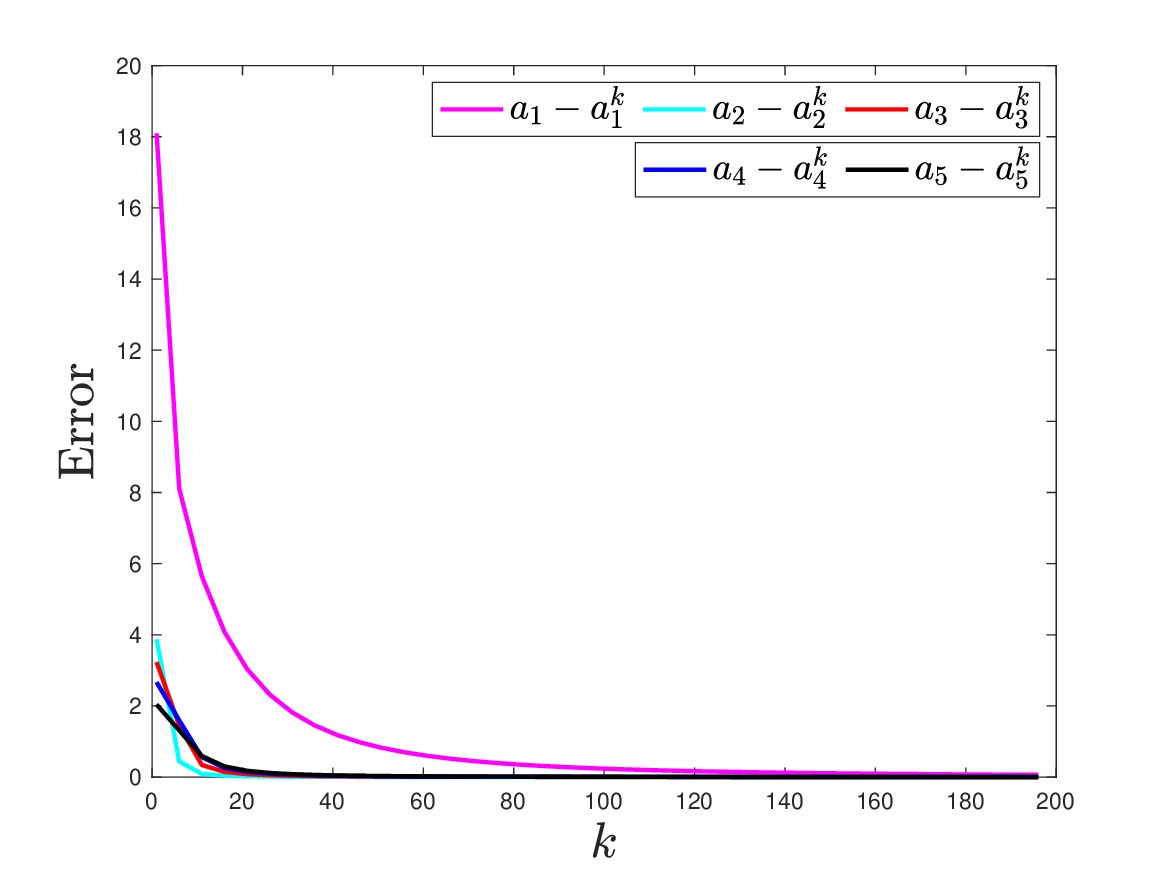}
&\includegraphics[width=0.33\textwidth]{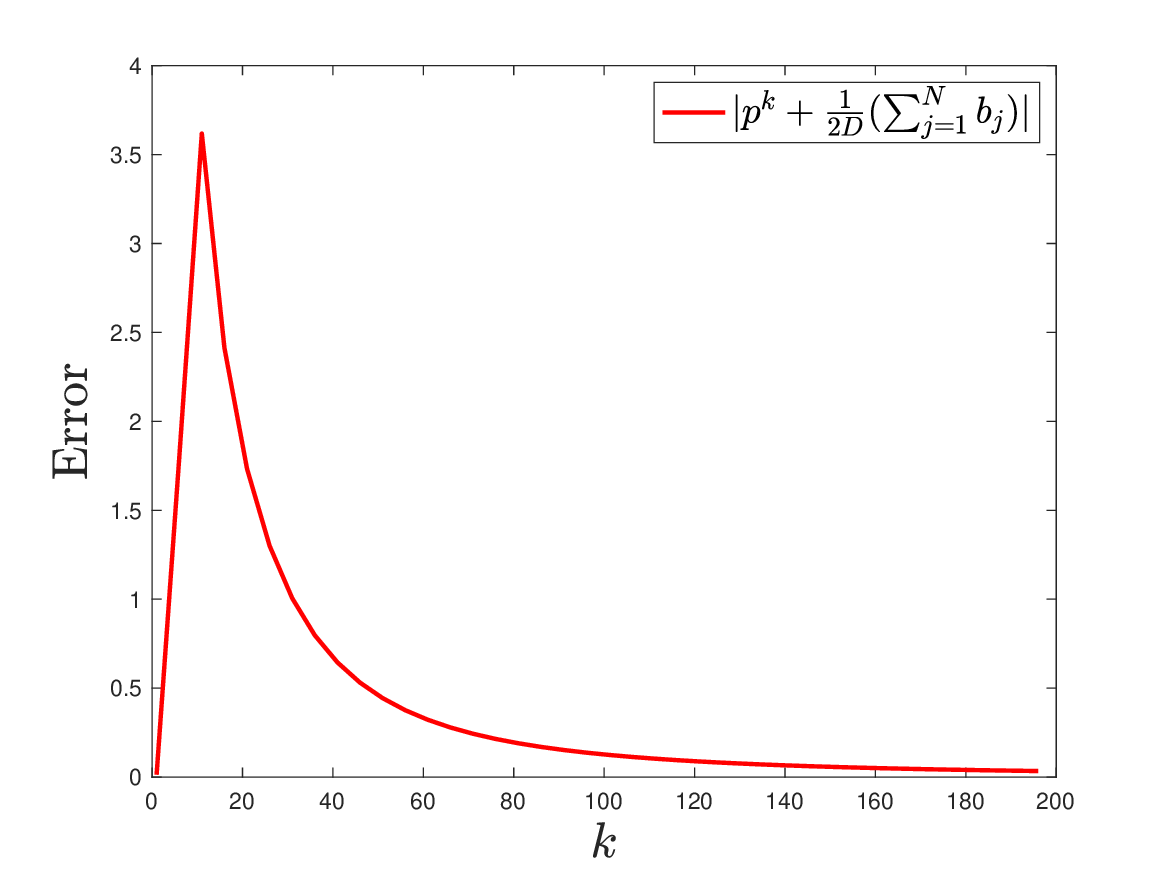}
&\includegraphics[width=0.33\textwidth]{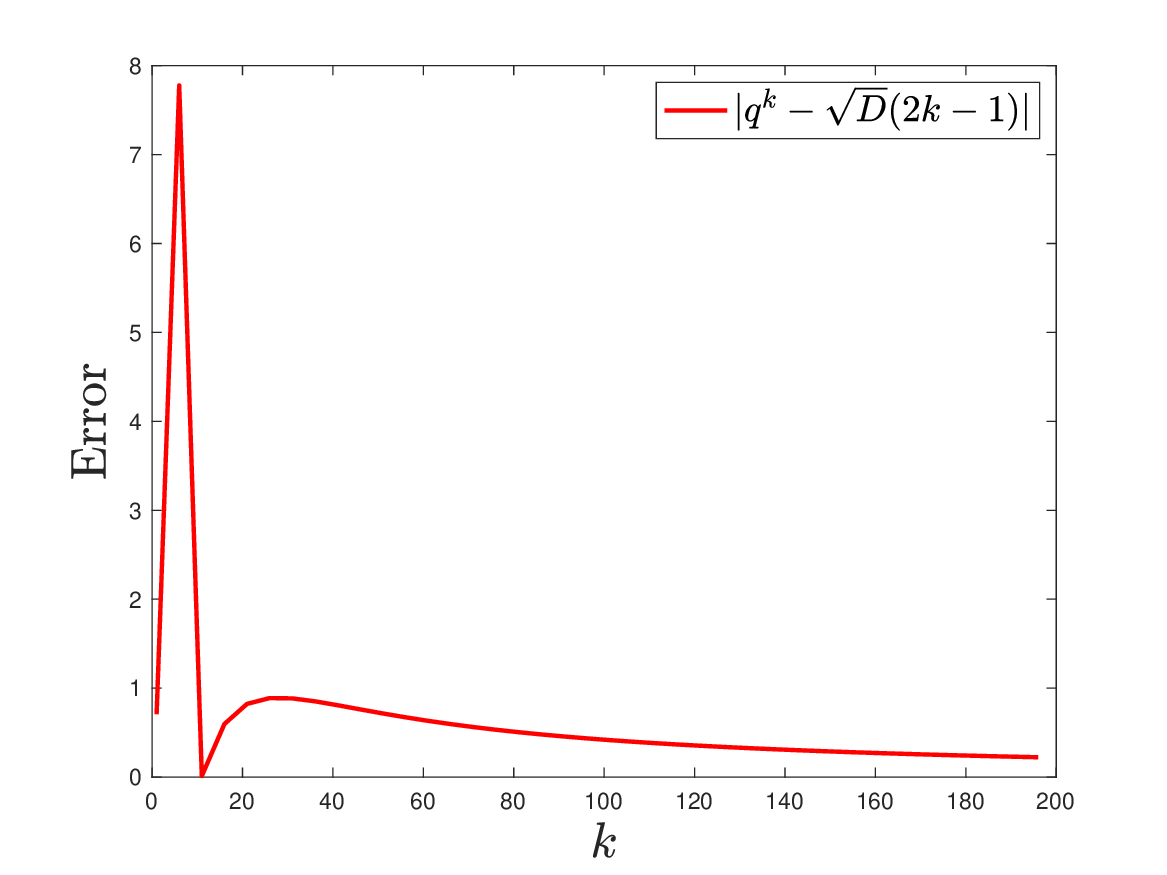}\\
(d) real roots for $D=1$ &(e) real part for $D=1$  &(f) imag. part for $D=1$ \\
\includegraphics[width=0.33\textwidth]{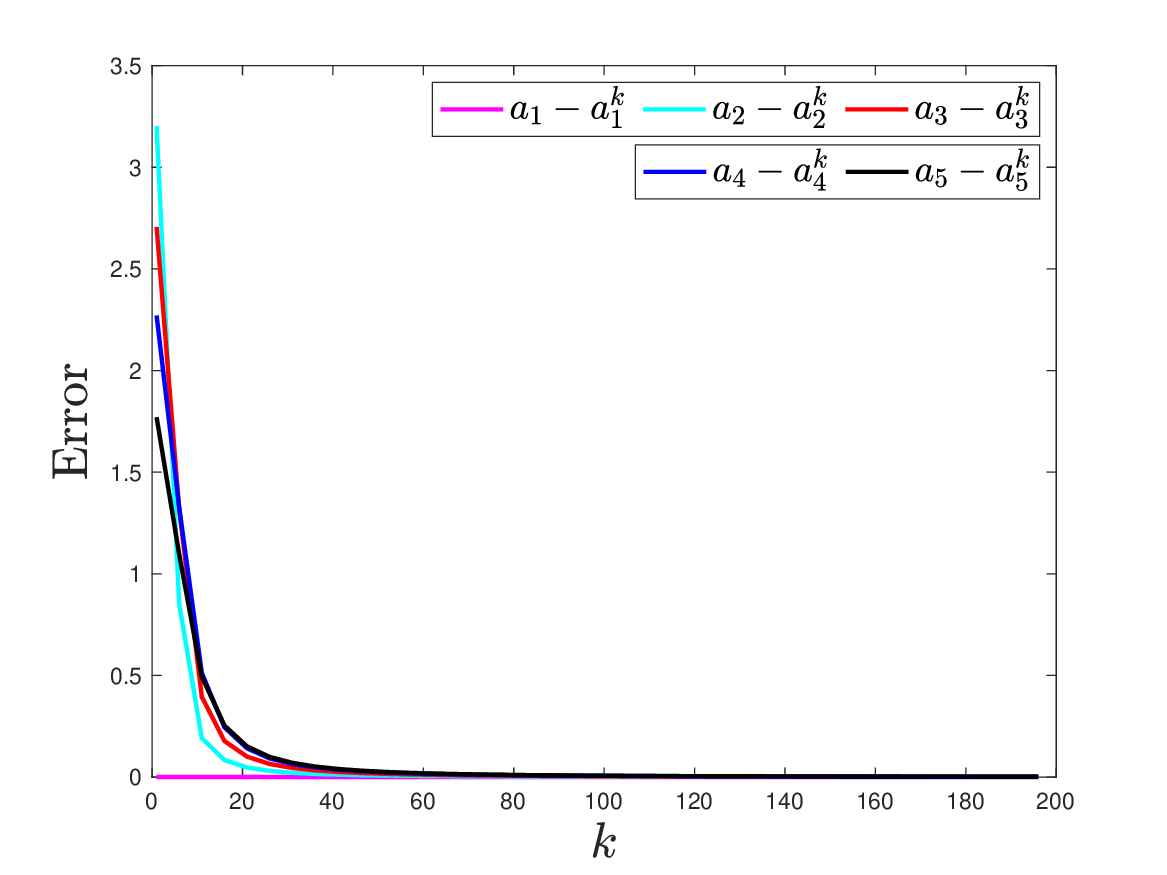}
&\includegraphics[width=0.33\textwidth]{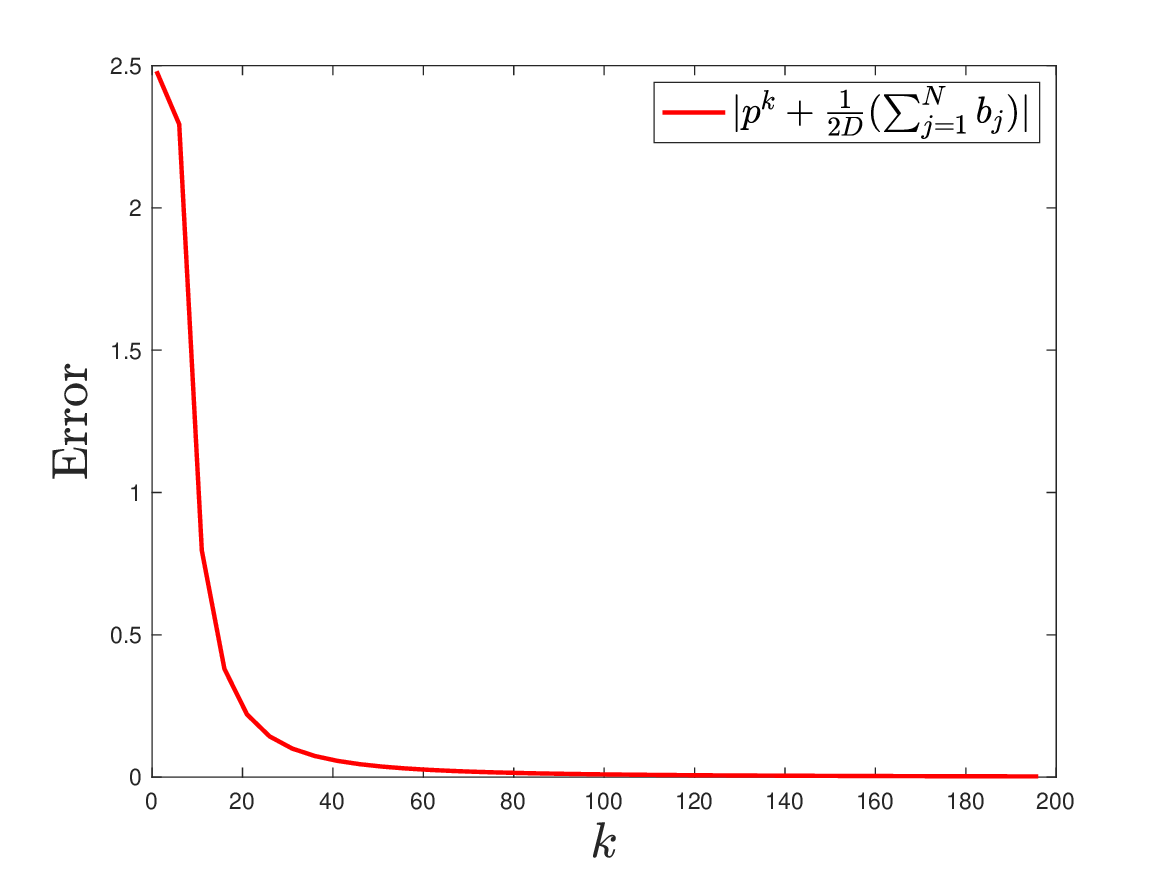}
&\includegraphics[width=0.33\textwidth]{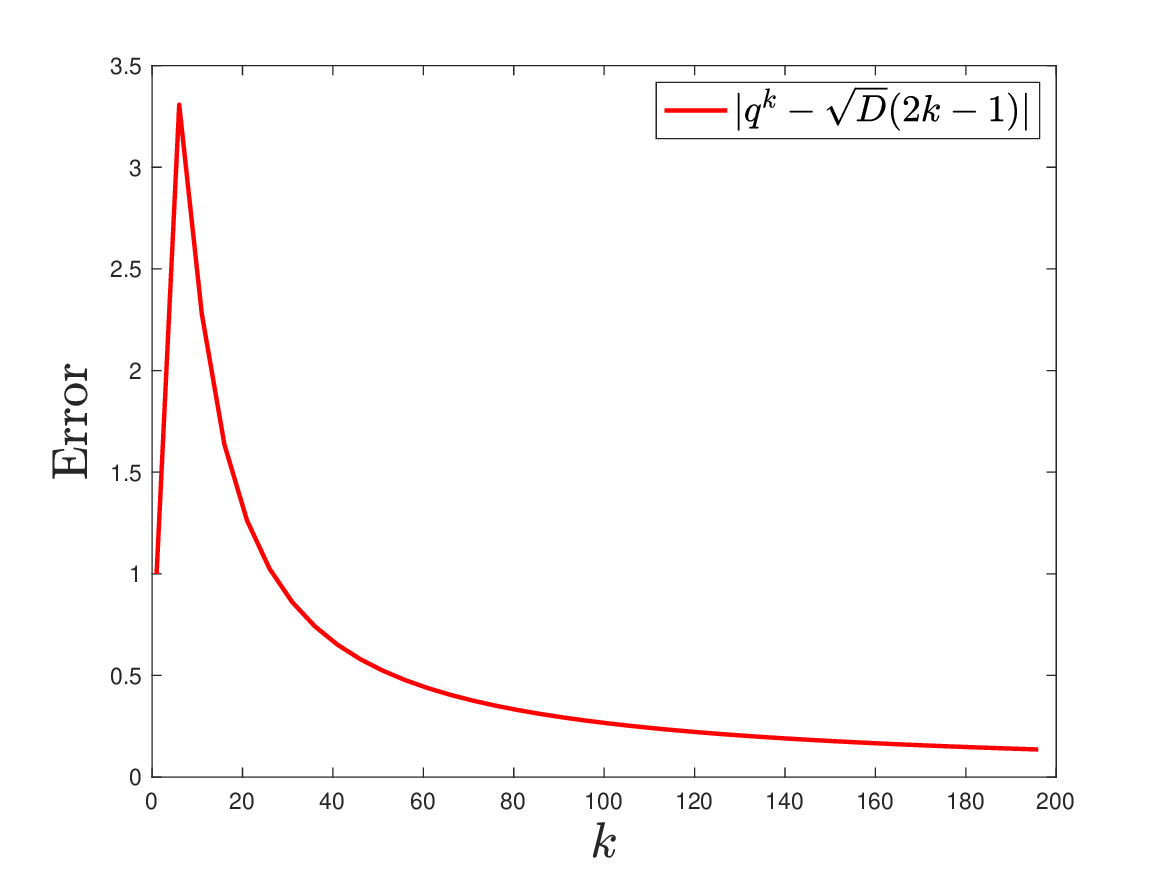}\\
(g) real roots for $D=5$ &(h) real part for $D=5$  &(i) imag. part for $D=5$  \\
\includegraphics[width=0.33\textwidth]{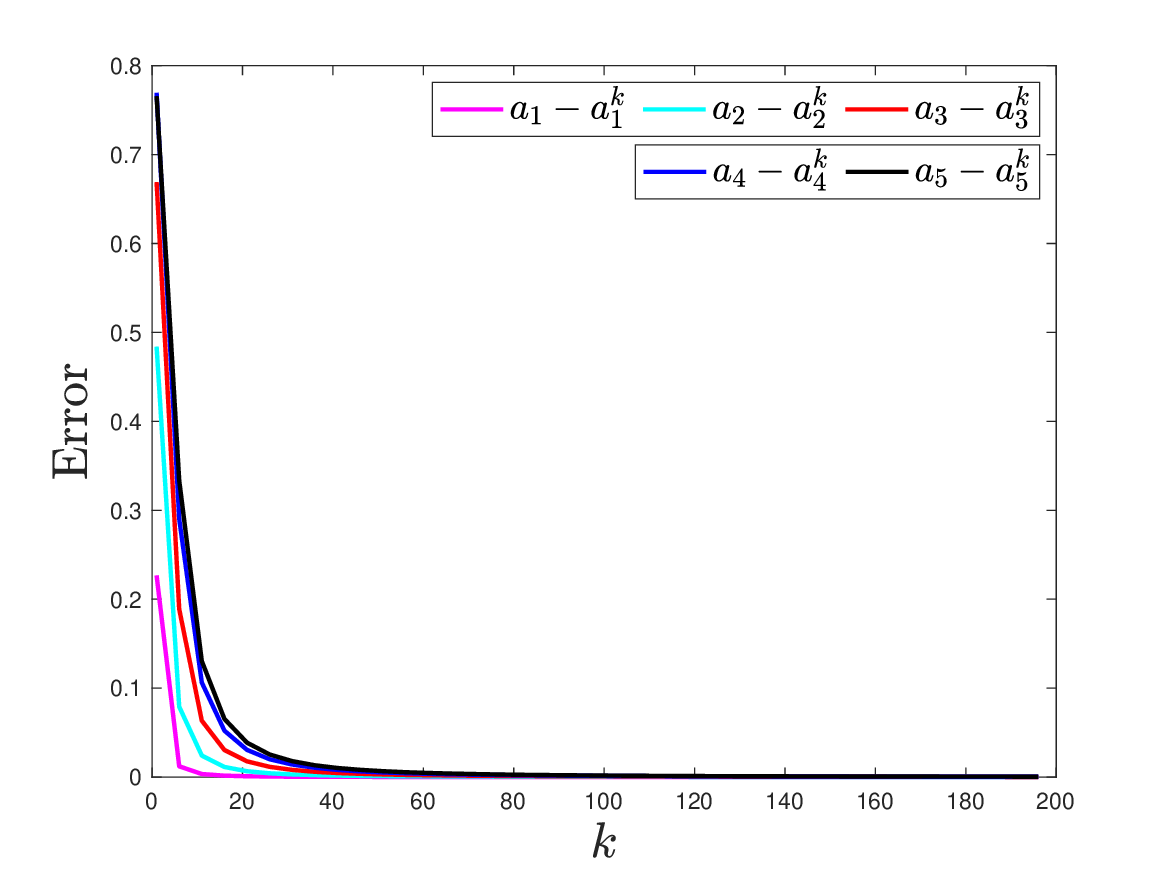}
&\includegraphics[width=0.33\textwidth]{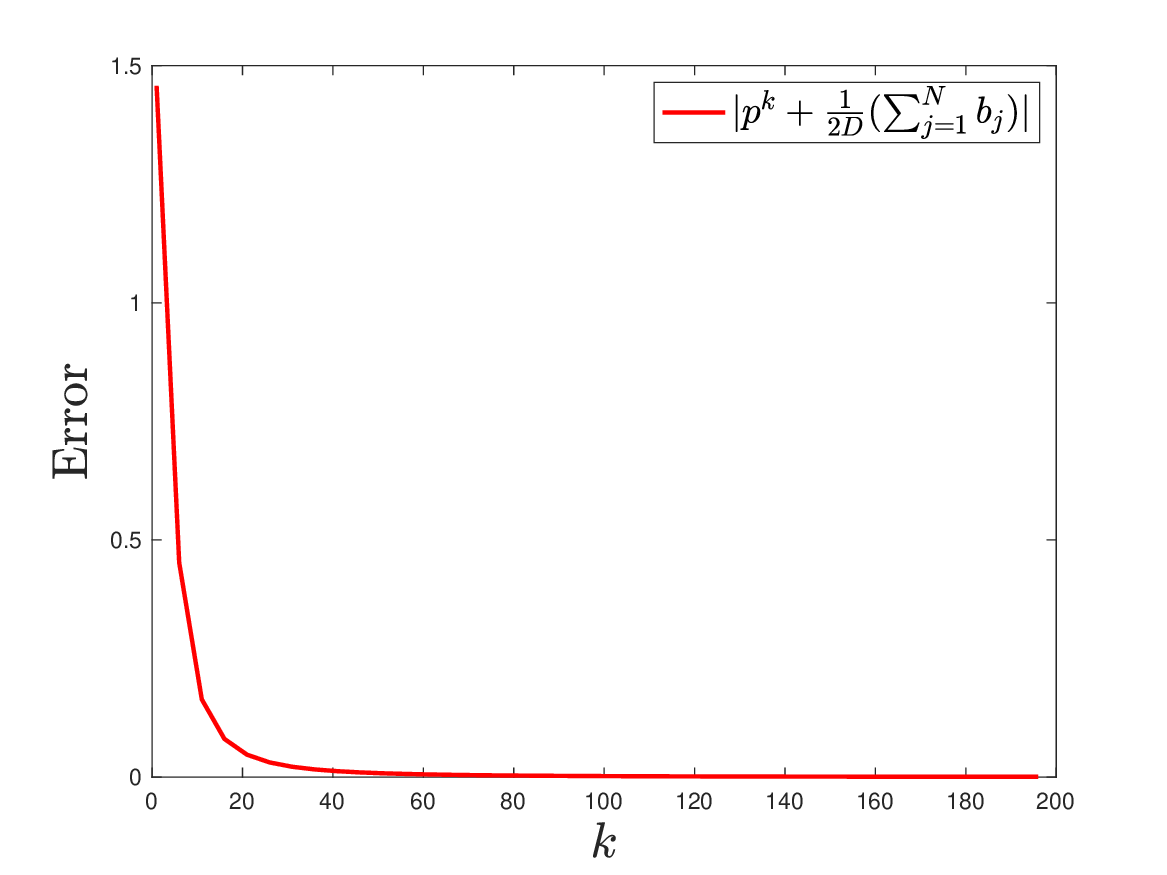}
&\includegraphics[width=0.33\textwidth]{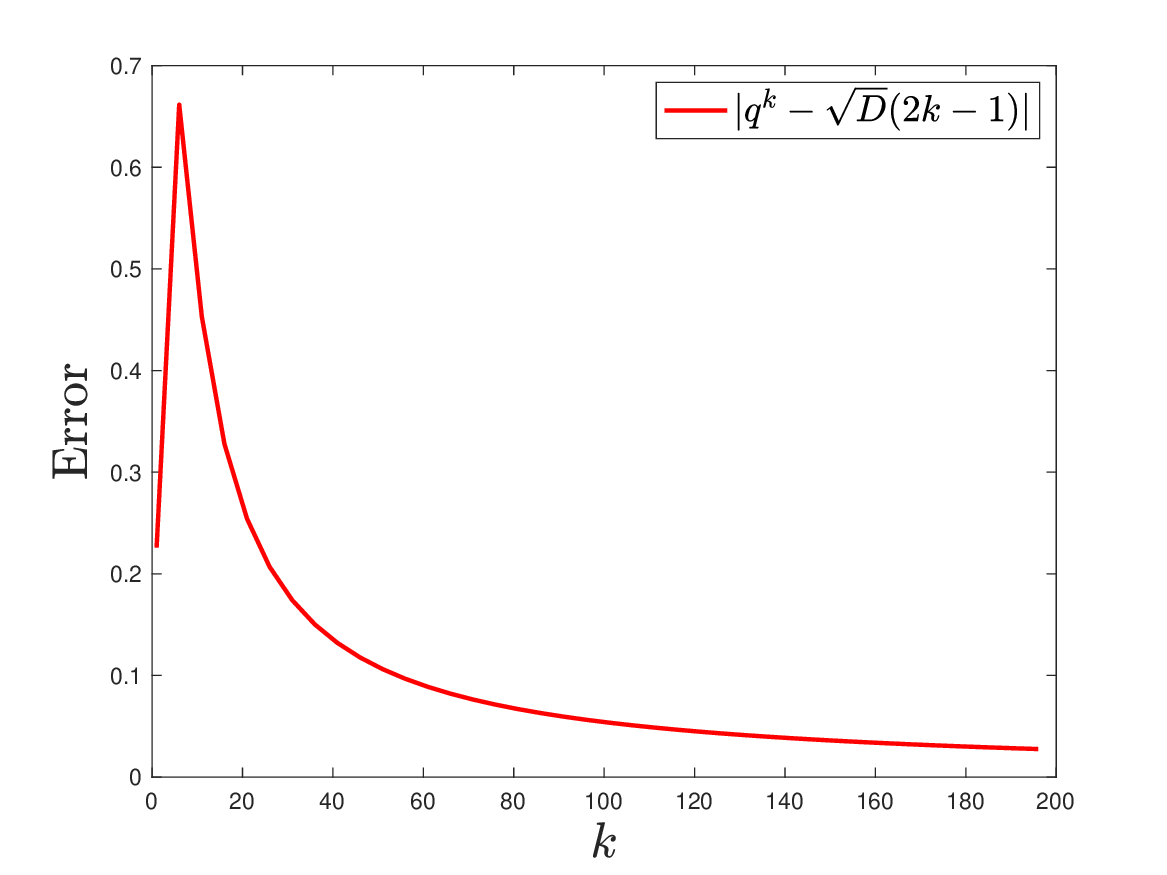}
\end{tabular}
\caption{Convergence of roots of $P_N^k$ for $N=5$ and $D=0.5,\,1,\,5$}
\label{fig_N5_Err}
\end{figure}

\begin{example}\label{exm_3}
Compared to the theoretical estimates of the convergence of roots given in \eqref{n5.16} and \eqref{n5.17}, we describe in Figure \ref{fig_N5_Err} the convergence of roots in more detail for the roots calculated in Example \ref{exm_1}.
\end{example}

%For both $D=1$ and $D=5$, the asymptotic behavior of the roots are similar. Taking $D=1$ as an example, we observed the asymptotic behavior of the real roots shown in Figures \ref{fig_N5} (a) and \ref{fig_N9} (a). Each real root of $P_N^k$, $a_j^k,\,j=1,\,\cdots,\,N$, converges to each real root of $P_N$, $a_j,\,j=1,\,\cdots,\,N$, as $k$ increases. In addition, the roots $a_j^k$ are located at $(-r_j,\,-r_{j-1})$ for $j=1,\,\cdots,\,N$. For the complex roots of $P_N^k$, there exists a pair of conjugate complex roots.  The real part and the imaginal part of complex roots satisfy the property given in Lemma \ref{lem5.3}.

In summary, we observed that the calculated roots are consistent with our theoretical results. More precisely, we first observed in Figures \ref{fig_N5g05_1_5}, \ref{fig_N9g05_1_5}, the consistency of calculated locations of the real roots $a_j$'s, $a_j^k$'s and complex roots $p_k\pm i q_k$ with those given in Lemmas \ref{lem4.1}, \ref{lem5.1}. Further, the consistency of calculated limits of $a_j^k,\,j=1,\,2,\,\cdots,\,N$ and $p_k\pm iq_k$ as $k$ increases numerically converge to $a_j,\,j=1,\,2,\,\cdots,\,N$ and $\sum_{j=1}^N\frac{b_j}{2D}\pm \sqrt{D}(2k-1)$, respectively. Then, in Figure \ref{fig_N5_Err}, we compared the numerical convergence rates of real roots and complex roots, and observed that the real roots have a faster convergence rate, which is consistent with Lemmas \ref{lem5.2}, \ref{lem5.3} saying that the convergence rate of the real root is $k^{-2}$, while for the complex roots $p_k\pm i q_k$, they are $k^{-2}$ and $k^{-1}$ for $p_k$ and $q_k$, respectively. 

\begin{comment} 
On the other hand, the convergence rates of both real roots and complex roots in the case of $D=0.5$ is significantly slower than that in the cases of $D=1$ and $D=5$.

Although the presented numerical results are limited, the correctness of the theoretical analysis presented in Section \ref{sec_asy} has been verified. That is, we numerically verify the number, the locations and the convergences of roots of $P_N^k$ defined by \eqref{char_poly}. 
\end{comment}

\section{Inverse Problem}
So far, we have obtained enough information on the structure of the characteristic roots of the characteristic polynomial \eqref{char_poly}, which is useful for data analysis of glassy state relaxation by looking at spectral data obtained by the mentioned DMA instruments. The next natural problem is to consider an inverse spectral problem. Based on the structure of the roots, we consider the following inverse problem in this section.

\medskip\noindent
{\bf Inverse problem:} Recover $D$ and the relaxation function given as Prony series \eqref{Prony} by knowing the eigenvalues of the reduced eigenvalue problems \eqref{new1.7} for several $k\in{\mathbb N}$. In other words, recover $D$ and the relaxation function by knowing several clusters of eigenvalues of \eqref{n1.7}. 

\medskip
An answer to this inverse problem is as follows.
\begin{theorem} By knowing two clusters of eigenvalues associated with $k=k_1,\,k_2\in{\mathbb N}$, we can recover $D$ and the relaxation function.
\end{theorem}
\pf.
The idea of proof is based on the well-known fact that for a monic polynomial, knowing all of its roots is equivalent to knowing the polynomial.
To begin further details of the proof, let $k_1, k_2\in{\mathbb N},\,k_1<k_2$ and assume that we know the clusters of eigenvalues associated to $k_1,\,k_2$ given as 
\begin{equation}\label{i7.1}
\begin{aligned}
a^{k_1}_1, a^{k_1}_2, a^{k_1}_3, \cdots,a^{k_1}_{N}, a^{k_1}_{N+1}=p^{k_1}+iq^{k_1},  a^{k_1}_{N+2}=p^{k_1}-iq^{k_1}
\end{aligned}
\end{equation}
and
\begin{equation}\label{i7.2}
\begin{aligned}
a^{k_2}_1, a^{k_2}_2, a^{k_2}_3, \cdots,a^{k_2}_{N}, a^{k_2}_{N+1}=p^{k_2}+iq^{k_2},  a^{k_2}_{N+2}=p^{k_2}-iq^{k_2},
\end{aligned}
\end{equation}
where $a^{k_1}_{N+1}$, $a^{k_1}_{N+2}$ or $a^{k_2}_{N+1}$, $a^{k_2}_{N+2}$ may be real.

Let us recall the characteristic polynomial for $P^k_N(\lambda)$ in \eqref{char_poly} given as
\begin{equation}\label{i7.3}
\begin{aligned}
P^k_N(\lambda)&=(D+\frac{\lambda^2}{(2k-1)^2})\Pi_{1\leq j\leq N}(\lambda+r_j) -  \sum_{i=1}^N b_i\Pi_{1\leq j\leq N,j\neq i}(\lambda+r_j)\\
&=\frac{1}{(2k-1)^2}\lambda^2\Pi_{1\leq j\leq N}(\lambda+r_j)+P_N(\lambda).
\end{aligned}
\end{equation}
Since this polynomial is of degree $N+2$ and its coefficient of $\lambda^{N+2}$ is $\frac{1}{(2k-1)^2}$, we have 
\begin{equation}\label{i7.4}
\begin{aligned}
P^{k_1}_N(\lambda)=\frac{1}{(2k_1-1)^2}\cdot\Pi_{1\leq j\leq N+2}(\lambda-a_j^{k_1})
\end{aligned}
\end{equation}
and
\begin{equation}\label{i7.5}
\begin{aligned}
P^{k_2}_N(\lambda)=\frac{1}{(2k_2-1)^2}\cdot\Pi_{1\leq j\leq N+2}(\lambda-a_j^{k_2})
\end{aligned}
\end{equation}
from \eqref{i7.1} and \eqref{i7.2}, respectively.
By putting $k=k_1$ and $k=k_2$ in \eqref{i7.3}, we have
\begin{equation}\label{i7.6}
\begin{aligned}
P^{k_1}_N(\lambda)-P^{k_2}_N(\lambda)&=\frac{1}{(2k_1-1)^2}\lambda^2\Pi_{1\leq j\leq N}(\lambda+r_j)-\frac{1}{(2k_2-1)^2}\lambda^2\Pi_{1\leq j\leq N}(\lambda+r_j)\\
&=[\frac{1}{(2k_1-1)^2}-\frac{1}{(2k_2-1)^2}]\cdot\lambda^2\Pi_{1\leq j\leq N}(\lambda+r_j).
\end{aligned}
\end{equation}
This implies
\begin{equation}\label{i7.7}
\begin{aligned}
\lambda^2\Pi_{1\leq j\leq N}(\lambda+r_j)=[\frac{1}{(2k_1-1)^2}-\frac{1}{(2k_2-1)^2}]^{-1}\cdot\big(P^{k_1}_N(\lambda)-P^{k_2}_N(\lambda)\big).
\end{aligned}
\end{equation}
Here, note that we know $k_1, k_2$ and all the respective roots of \eqref{i7.4} and \eqref{i7.5}, the polynomial on the right-hand side
of \eqref{i7.7} is known due to the mentioned well-known fact. Hence we know $r_j,\,1\leq j\leq N$.

Now, recall \eqref{n4.8} given as
\begin{equation*}
\begin{aligned}
P^k_N(-r_i)= -   b_i\Pi_{1\leq j\leq N,j\neq i}(-r_i+r_j)
\end{aligned}
\end{equation*}
for $1\leq i\leq N$.
Using this for $k=k_2$, we can recover $b_i$
for $1\leq i\leq N$ because we know $P_N^{k_2}(\lambda)$ and $r_i,\,1\leq i\leq N$. We can further recover $D$ using \eqref{i7.3} with $k=k_2$ because we already know $b_i,\,r_i,\,1\le i\le N$. Thus, we have recovered $D$, $r_i$ and $b_i$ for $1\leq i\leq N$ by knowing the eigenvalues of $P^{k_1}_N(\lambda)$ and $P^{k_2}_N(\lambda)$ with $k_1, k_2 \in \mathbb{N}$, $k_1<k_2$. 
\eproof

\section{Discussion}
In this paper, we provided  rigorous insight in clusters of eigenvalues associated with an extended Burgers model for viscoelastic relaxation and the structure of the characteristic roots of the characteristic polynomial for the reduced eigenvalue problem. The very interesting results that have emerged from these are as follows. We emphasize the appearance of quasi-static modes, which are purely damping modes for the mentioned insight. Since this mode appears in the free oscillation of the Earth, it suggests that the method of this paper may be useful for analyzing the free oscillation of the Earth (see the next paragraph for further details). As for the mentioned structure, we emphasize the success in showing that we can recover $D$ and the relaxation function by knowing two clusters of eigenvalues. This is very important, because it can initiate a robust data analysis for glassy state relaxation, and also the recovered $D,\,b_i\text{'s},\,r_i\text{'s}$ are translated into the stretched  exponential function approximately by using the mentioned approximate correspondence \cite{MM}. 
\begin{comment}
The details will be given in a forthcoming paper.    
\end{comment}

\begin{comment}
Lemmas \ref{lem4.1}, \ref{lem4.2}, \ref{lem5.1}, \ref{lem5.2}, \ref{lem5.3} and their proofs not only provide this insight but also indicate the possibility of identifying $D$ and the relaxation from eigenvalues. In fact, using two clusters of eigenvalues associated with $k=k_1,\,k_2\in{\mathbb N}$, we can recover $D$ and the relaxation. 
As the eigenvalues of $A_{N+2}^k$ are simple roots, we can approximately recover $D,\,b_i\text{'s},\,r_i\text{'s}$. A brief description is as follows. From the asymptotic behavior of $q^k/(2k-1)$, we can assume that $D$ is approximately known. Generically, each $a_j^k\text{'s}$ and $b_j\text{'s}$ can be expressed in terms of $(b_i,r_i)\text{'s}$ and $(a_j^k,r_j)\text{'s}$, respectively. Combining these two expressions, we expect to have implicit equations for $(a_j^k,r_j)\text{'s}$ solvable with respect to $r_j\text{'s}$.
\end{comment}

We have focused on the asymptotic behavior of large eigenvalues. It is also feasible to analyze the behavior for low-lying eigenvalues. This is interesting, for example, in the study for free oscillations of the Earth, while it has been proposed that the EBM is an appropriate model for viscoelastic relaxation (see \cite{YP}). The rigorous study of quasistatic (toroidal or spheroidal) modes in this context has not yet been developed.

\subsection*{Acknowledgments}

We appreciate a useful comment by Kazumi Tanuma to improve our paper. As for the funding, the second author was supported by the Simons Foundation under the MATH + X program, the National Science Foundation under grant DMS-2108175, and the corporate members of the Geo-Mathematical Imaging Group at Rice University. The third author was supported by NSFC-RGC Joint Research Grant No. 12161160314 and the Natural Science Foundation Innovation Research Team Project of Guangxi (Grant No.2025GXNSFGA069001). The fourth author was partially supported by the Ministry of Science and Technology of Taiwan (Grant No. NSTC 111-2115-M-006-013-MY3).
The fifth author was partially supported by JSPS KAKENHI (Grant No. JP22K03366, JP25K07076).

\end{document}